\documentclass[12pt,reqno]{amsart}
\usepackage{amsmath}
\usepackage{amssymb}
\usepackage{verbatim}
\usepackage{graphicx}
\usepackage[font=footnotesize]{caption}
\usepackage{subcaption}
\usepackage{epstopdf}
\usepackage{multirow}
\usepackage{url}
\usepackage{makecell}
\usepackage{tabularx}
\usepackage{longtable}
\usepackage{graphicx}
\usepackage{tikz}
\usepackage{pgfplots}
\usepgfplotslibrary{polar}
\usepackage{xcolor}
\usepackage{pgfplots}
\usepackage{pgfplotstable}
\allowdisplaybreaks

\usepackage{extarrows}
\usepackage{bm}
\usepackage[toc]{appendix}
\usepackage[capitalise]{cleveref}
\usepackage{color}
\allowdisplaybreaks

\usepackage[top=0.5in,bottom=0.5in,left=0.7in,right=0.7in]{geometry}

\usepackage{bm}

\newtheorem{theorem}{Theorem}[section]
\newtheorem{remark}{Remark}[section]
\newtheorem{example}{Example}[section]
\numberwithin{equation}{section}
\numberwithin{lemma}{section}

\numberwithin{equation}{section}

\newcommand{\N}{\mathbb{N}}    
\newcommand{\R}{\mathbb{R}}    

\newcommand{\bo}{\mathcal{O}}

\newcommand{\be}{ \begin{equation} }
	\newcommand{\ee}{ \end{equation} }

\newcommand{\nab}{\nabla}

\begin{document}

		\title[]{Efficient and simple fourth-order compact finite difference methods  for  convection-diffusion-reaction equations on arbitrary curved domains}
		
		\author{Qiwei Feng, Bin Han, and Peter Minev}
		\address{Mathematics Department, King Fahd University of Petroleum and Minerals, Dhahran, 31261, Saudi Arabia.}
		\email{ {\tt qiwei.feng@kfupm.edu.sa,qfeng@ualberta.ca}}
		
	\address{Department of Mathematical and Statistical Sciences,
		University of Alberta, Edmonton, Alberta, Canada  T6G 2N8.}
	\email{ {\tt bhan@ualberta.ca}
		\quad {\tt minev@ualberta.ca}}

		\makeatletter \@addtoreset{equation}{section} \makeatother

		\begin{abstract}	

	In this paper, we discuss the 2D convection-diffusion-reaction equation with variable smooth coefficients and the Dirichlet boundary condition on a complicated, thin, and curved domain. We propose the fourth-order compact FDM at every grid point with the uniform Cartesian mesh. For the regular stencil center, we utilize the fourth-order compact 9-point FDM to approximate the solution. According to the preliminary analysis, we use
	 vertical and horizontal transformations to derive fourth-order compact FDMs in 10 cases for all irregular  stencil centers. To obtain the left-hand side of the stencil of the fourth-order FDM in each case, we only need to solve an at most $6 \times 24$ linear system which is presented with the explicit formula. The right-hand side of the FDM is constructed in explicit expression for any irregular  stencil centers too. To achieve the fourth-order consistency, up to second-order partial derivatives of convection, diffusion, reaction, and source terms are used for the FDM at the regular  stencil center, and the FDM at an irregular  stencil center
	 only requires first-order partial derivatives of convection, diffusion, reaction, and source terms, and up to third-order derivatives of the Dirichlet boundary function and the parametric expression of the boundary curve. We test  challenging domains with 100-leaf, high-curvature, high-frequency, sharply varying, and nearly overlapping boundary curves, the proposed FDM produces the high accuracy and the stable fourth-order convergence rate in $l_2$ and $l_{\infty}$ norms. From error plots, we observe that errors are uniformly distributed on the closure of domains which confirms the robustness and stability.
	  All examples   numerically indicate that  the accuracy and convergence rate seem to be independent of the distance between the boundary and the irregular  stencil center when the mesh size $h$ is reasonably small.
 All stencils of our FDMs have a simple desired structure by only keeping grid points inside $\Omega$ in the standard compact 9-point stencil for both regular stencils and boundary stencils, but without assuming any information outside the domain $\Omega$.
		\end{abstract}	
\maketitle
\vspace{-0.6cm}

\pagenumbering{arabic}

	\section{Introduction}

Partial differential equations (PDEs) on the complex irregular domain with the arbitrary geometry
 arise in many scientific and engineering applications. Typical examples include the fluid flow over the  arbitrary domain, the heat transfer in the heterogeneous material, the free boundary problem, the parabolic moving boundary problem, and the diffusion process around the boundary with the irregular geometry. The complicated geometry of the irregular domain significantly increases the computational cost to obtain the accurate solution, because the sharply varying boundary requires very fine mesh. So in this paper, we derive the high-order finite difference method (FDM) to handle the arbitrary curved boundary. Since this paper focuses on FDMs, we mainly review the following literature of deriving efficient FDMs for various crucial PDEs over irregular domains, other numerical methods can be found in the references therein.

Arias  et al. in \cite{Arias2018} proposed the conservative method with the second-order  accuracy for  $-\nab\cdot (\alpha \nab u)=\phi$ with the constant diffusion term $\alpha$ and the Robin boundary condition $\rho u+ \alpha \frac{\partial u}{\partial \vec{n}}=g$ on the 2D irregular domain. Artzi et al. in \cite{Artzi2009} proved the  second-order convergence rate of a compact FDM for the biharmonic equation on the 2D planar irregular domain.
Baeza  et al. in \cite{Baeza2016} constructed the high-order extrapolation technique of the weighted essentially non-oscillatory (WENO) FDM at the boundary of the 2D complex domain for the hyperbolic conservation law.  The second-order unified FDM to solve the Poisson equation and the heat transfer problem  with Dirichlet, Neumann, and Robin  boundary conditions on the 2D
irregular  domain was discussed by Chai et al. in \cite{Chai2021}.
The fast  second-order FDM for the biharmonic equation on the 2D irregular domain was derived by Chen  et al. in \cite{Chen2008}.
 The second-order FDM for variable coefficient Poisson and heat equations with the Dirichlet boundary condition using the non-graded adaptive  grid on 2D and 3D irregular domains was established by Chen  et al. in \cite{Chen2007}. Clain	et al.  proposed the arbitrary high-order FDM for the linear convection-diffusion equation \cite{Clain2021} and the nonlinear convection-diffusion-reaction equation \cite{Clain2024} with Dirichlet, Neumann, linear or nonlinear Robin conditions on the arbitrary smooth 2D domain.
Ding et al.	applied the stencil-adaptive FDM to solve the 2D
unsteady incompressible viscous flow with the curved boundary in \cite{Ding2007}.
The fourth-order embedded FDM for the 2D initial boundary value problem with the moving boundary and the Dirichlet boundary condition was built by Ditkowski et al. in \cite{Ditkowski2009}. Fidalgo	et al. in \cite{Fidalgo2020} proposed the fifth-order WENO FDM with the ghost cell method  to approximate the solution of the conservation law on the curved boundary domain. The fourth-order FDM for  Laplace and heat equations with the Dirichlet
boundary condition on the 2D arbitrary domain was designed by Gibou et al. in \cite{Gibou2005}.
The second-order symmetric FDM for variable coefficient  Poisson and heat equations with the Dirichlet
boundary condition on 2D and 3D irregular domains was presented by Gibou et al. in \cite{Gibou2002}. Han  et al. in \cite{HanSim2025} introduced the sixth-order compact FDM for the variable coefficient  Poisson equation with the Dirichlet
boundary condition on the 2D curved domain, and  provided the corresponding convergence analysis.
Helgad\'{o}ttir et al. in \cite{Helgadottir2011} deduced the  second-order FDM for the nonlinear Poisson-Boltzmann equation with Neumann or Robin
boundary condition on 2D and 3D irregular domains with the arbitrary geometry
 using the non-graded adaptive grid, and the corresponding linear system in \cite{Helgadottir2011} to compute the solution is symmetric positive definite if the uniform grid is applied. By employing 1 to 4 points at the boundary,  Ito  et al. in \cite{Ito2005} developed
 the fourth-order compact FDM or immersed interface method (IIM) for the diffusion-reaction equation on the irregular domain with the Robin boundary condition, and Li  et al. in \cite{LiIto06} derived the fourth-order compact FDM/IIM for Poisson and heat equations on the 2D irregular domain with the Dirichlet boundary condition (necessary steps are provided to deal with the 3D irregular domain). The embedded FDM for the Poisson equation  with the Dirichlet
boundary condition on the 2D irregular domain was proposed by Jomaa  et al. in \cite{Jomaa2005}. Lakner et al. in \cite{Lakner2008}  established the  efficient symbolic-numerical procedure of the FDM for Laplace and heat equations with Dirichlet and Neumann boundary conditions
over the 2D curvilinear domain. Li  et al. in \cite{Li2023,Li2025}  formulated the fourth-order in space and the second-order in
time augmented matched interface and
boundary (AMIB) method with the Fast Sine or Fast Fourier Transform acceleration for elliptic and/or parabolic problems over 2D and/or 3D irregular domains with   Dirichlet, Neumann, and Robin boundary conditions. Liszka et al. in \cite{Liszka1980} applied the FDM with the arbitrary irregular grid to solve the applied mechanics problem. Ng  et al. in \cite{Ng2009} proposed the second-order ghost fluid method for the variable coefficient Poisson equation with the Dirichlet boundary condition on 2D and 3D irregular domains. Olsson in	 \cite{Olsson1994} derived the FDM with the third-order accuracy at the boundary and the sixth-order accuracy in the interior domain for parabolic and symmetric hyperbolic systems over the non-smooth domain. Pan  et al. in \cite{Pan2021}  developed an augmented third-order compact FDM/IIM
for the Poisson/Helmholtz equation with Dirichlet, Neumann, and Robin  boundary conditions on the 2D irregular domain.
The fourth-order AMIB method with the FFT acceleration for $\Delta u+\kappa u=\phi$ with the boundary condition  $\rho u+ \zeta \frac{\partial u}{\partial \vec{n}}=g$ on  2D and 3D irregular domains was established by Ren et al. in \cite{Ren2022}. The sharp second-order finite volume and difference method with the adaptive grid for the linear elasticity equation with Dirichlet and Neumann boundary conditions on 2D and 3D irregular domains was constructed by Theillard et al. in \cite{Theillard2013}. Yoon et al. in \cite{Yoon2015}  derived the second-order convergence proof of the FDM proposed by Gibou et al. for the Poisson equation with the Dirichlet boundary condition on the 2D irregular domain. The second-order augmented  IIM  for the Helmholtz/Poisson equation with the Dirichlet boundary condition on the 2D irregular domain in complex space was built by Zhang et al. in \cite{Zhang2016}. For the convergence proof of the high-order compact FDM over rectangular and cubic domains, Berikelashvili et al. in \cite{Berikelashvili2007} provided the convergence proof of the fourth-order compact FDM for the convection-diffusion equation
with the constant coefficient and the Dirichlet boundary condition in the unit cube, and
Shi   et al. in \cite{Shi2021} proved the stability and the optimal fourth-order convergence rate of the compact block-centered FDM for elliptic and parabolic problems with the variable coefficient on the rectangular domain.

\begin{figure}[htbp]
	\centering
	 \begin{subfigure}[b]{0.32\textwidth}
		 \includegraphics[width=5.5cm,height=5.5cm]{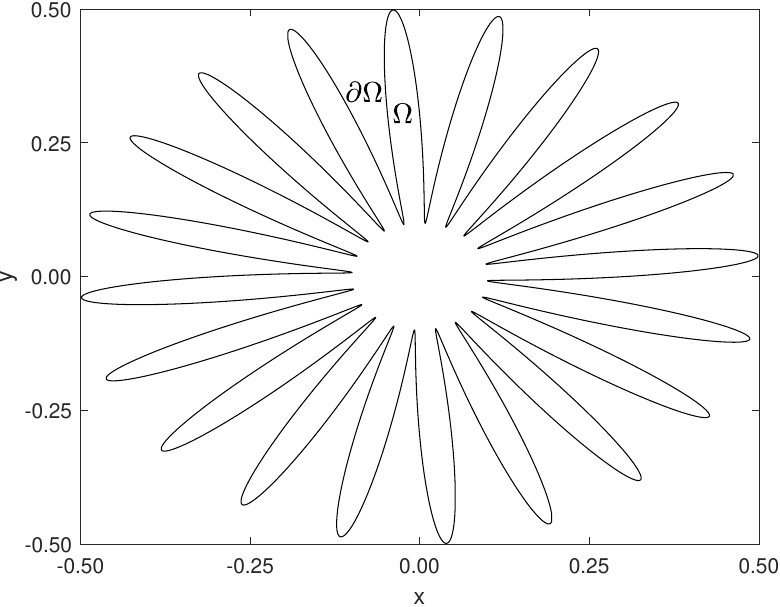}
	\end{subfigure}
	 \begin{subfigure}[b]{0.32\textwidth}
		 \includegraphics[width=5.5cm,height=5.5cm]{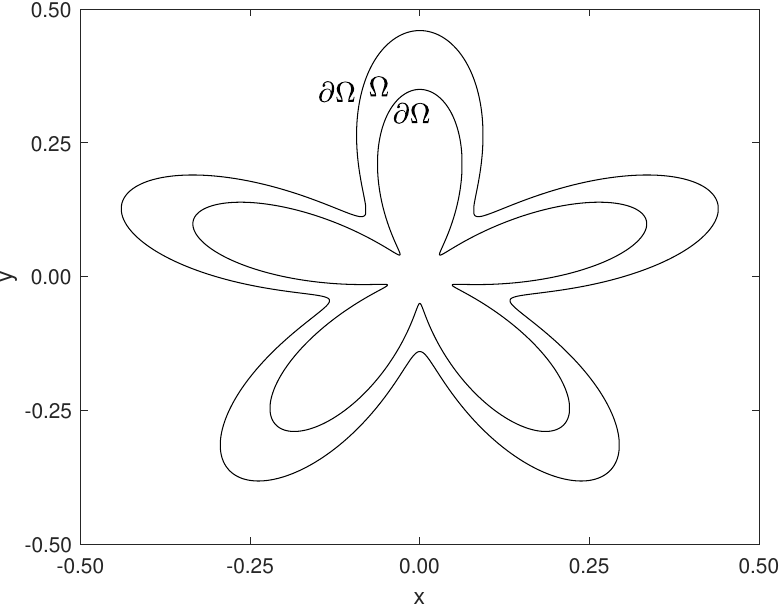}
	\end{subfigure}
	 \begin{subfigure}[b]{0.32\textwidth}
		 \includegraphics[width=5.5cm,height=5.5cm]{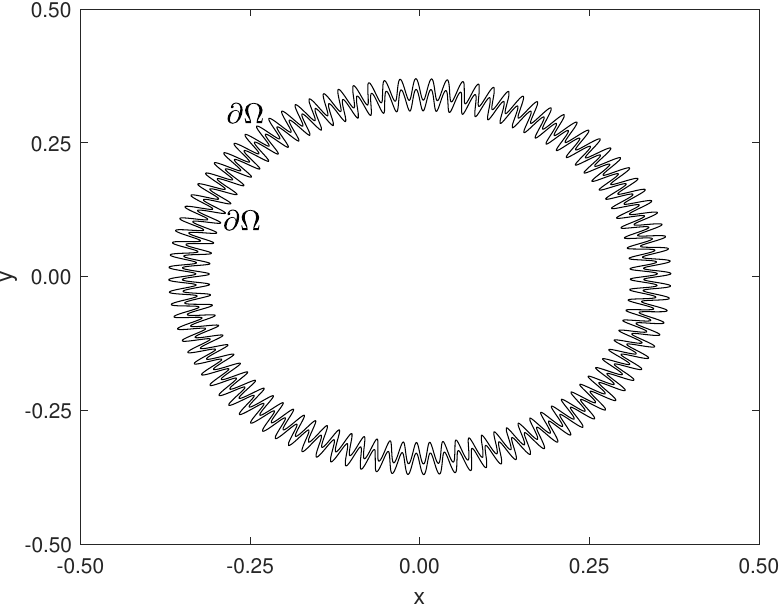}
	\end{subfigure}
	\caption
	{Three examples of $\Omega$ with complicated boundary curves $\partial \Omega$   for the model problem \eqref{model:problem:2D}: $\Omega$ is the region enclosed by a 20-leaf boundary curve (left),  $\Omega$ is the region enclosed by two 5-leaf boundary curves (middle), and $\Omega$ is the thin region enclosed by two 100-leaf boundary curves (right).}
	\label{illustra:fig}
\end{figure}	
	In this paper, we derive the fourth-order compact finite difference method (FDM) to solve the 2D convection-diffusion-reaction equation with the Dirichlet boundary condition on an arbitrary  irregular domain $\Omega\subset \R^2$ with the smooth boundary curve $\partial \Omega$ as follows:
\be\label{model:problem:2D}
\begin{cases}
	-\nab\cdot (\alpha \nab u) +    \vec{\beta}\cdot \nabla u +   \kappa u = \phi, \qquad \qquad  &  (x,y)\in \Omega,\\
	u =g,  & (x,y)\in \partial \Omega,
\end{cases}
\ee
where $\alpha>0$, $\vec{\beta}=(\beta_1,\beta_2)$, $\kappa$, and $\phi$ are variable functions in $\overline{\Omega}$.
See \cref{illustra:fig} for illustrations of $\Omega$ that we discuss in this paper to  construct the fourth-order compact FDM to solve \eqref{model:problem:2D}. To derive the fourth-order compact FDM, we assume that
\begin{itemize}
	\item the partial derivatives of  $u$ up to the fourth order are uniformly continuous on $\overline{\Omega}$,
	\item the partial derivatives of  $\alpha$, $\vec{\beta}$,  $\kappa$, and $\phi$ up to the second order are uniformly continuous on $\overline{\Omega}$,
	\item the derivatives of the Dirichlet boundary function $g$ and the parametric expression of the boundary curve $\partial \Omega$  up to the third order are uniformly continuous on $\partial \Omega$.
\end{itemize}
	Let
\be\label{notation:abdf}
a:=\frac{\alpha_x-\beta_1}{\alpha}, \qquad b:=\frac{\alpha_y-\beta_2}{\alpha}, \qquad d:=-   \frac{\kappa}{\alpha}, \qquad f:=-\frac{\phi }{\alpha}.
\ee
Then \eqref{model:problem:2D} implies
\be\label{2D:linear:simple:eq}
\Delta u  + a u_x +  bu_y +du =f.
\ee

\begin{figure}[htbp]
		\centering
		 \begin{subfigure}[b]{0.4\textwidth}
			 \includegraphics[width=7cm,height=7cm]{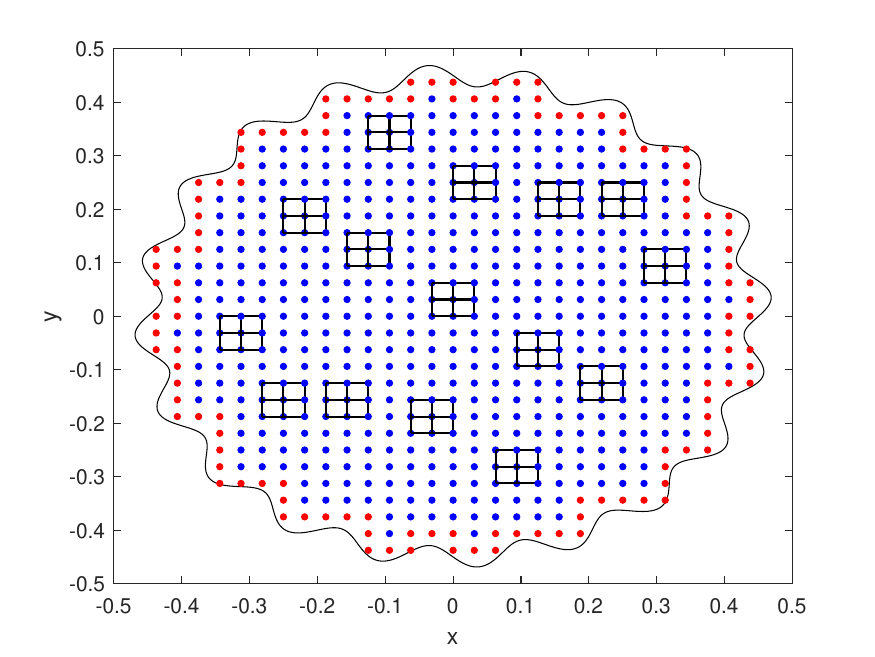}
		\end{subfigure}
		 \begin{subfigure}[b]{0.4\textwidth}
			 \includegraphics[width=7cm,height=7cm]{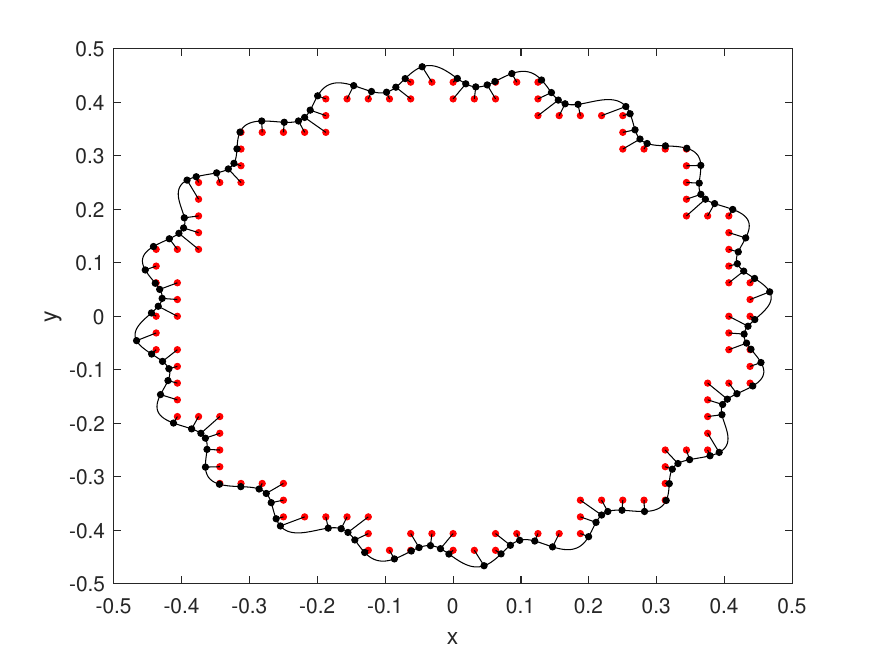}
		\end{subfigure}
		\caption
		{Left: Regular (blue) and irregular (red)  stencil centers. Right: Orthogonal projections (black points) of the stencil centers (red points) of  boundary irregular stencils onto the boundary curve $\partial \Omega$ (the black curve).}
			 \label{nine:point:orthog:fig}
	\end{figure}

For simplification, we assume that the original point $(0,0)\in \Omega$ to describe the uniform Cartesian grid used in this paper for $\Omega$:
\be \label{xiyj:2D:space}
(x_i,y_j) \cap \overline{\Omega} \quad \text{with} \quad (x_i,y_j):=(i h,jh), \quad i,j=0,\pm 1, \pm 2 \ldots,  \quad \text{and} \quad h:=1/N, \quad N\in \N.
\ee
We also define that
\[
u_{i,j}=u(x_i,y_j), \qquad (u_h)_{i,j}=\text{the value of } u_h \text{ that computed at the grid point } (x_i,y_j),
\]
where $u_h$ is the numerical solution that  approximated by our proposed fourth-order compact FDM for \eqref{model:problem:2D}.
Furthermore, we say that the compact FDM utilizes the regular stencil at the center point $(x_i,y_j)$ if $\{ (x_{i+\mu},y_{j+\nu}) : -1\le \mu,\nu\le 1\} \in \overline{\Omega}$ with  $(x_i,y_j)$ denoting the regular  stencil center (see blue points in the left panel of \cref{nine:point:orthog:fig} for illustrations), otherwise the compact FDM utilizes the irregular stencil with $(x_i,y_j)$ denoting the irregular  stencil center (see red points in the left panel of \cref{nine:point:orthog:fig} for illustrations). 
\begin{figure}[htbp]
	\centering
	\begin{subfigure}[b]{0.3\textwidth}
		 \includegraphics[width=6cm,height=6cm]{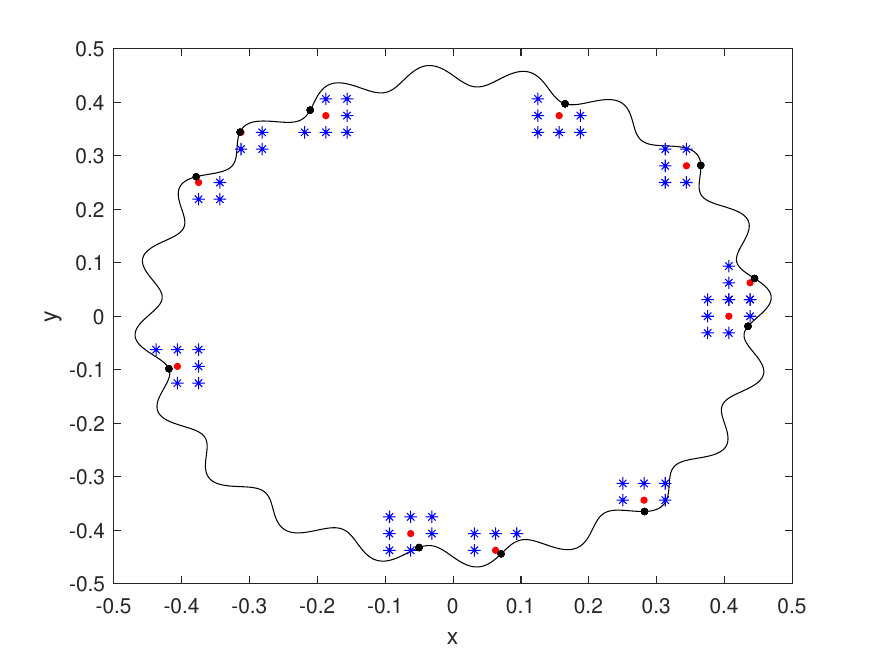}
	\end{subfigure}
	\begin{subfigure}[b]{0.3\textwidth}
		 \includegraphics[width=6cm,height=6cm]{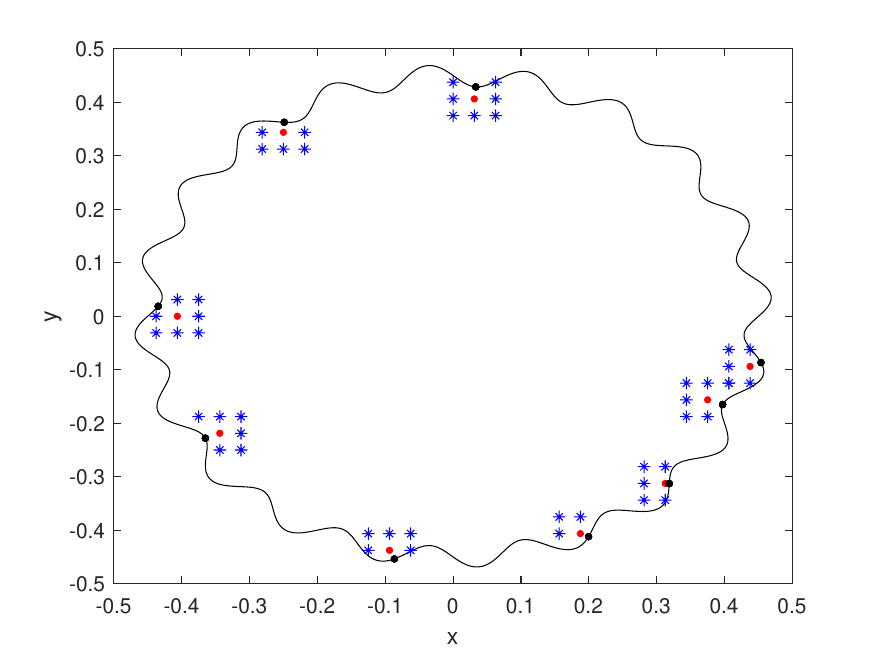}
	\end{subfigure}
		 \begin{subfigure}[b]{0.3\textwidth}
		 \includegraphics[width=6cm,height=6cm]{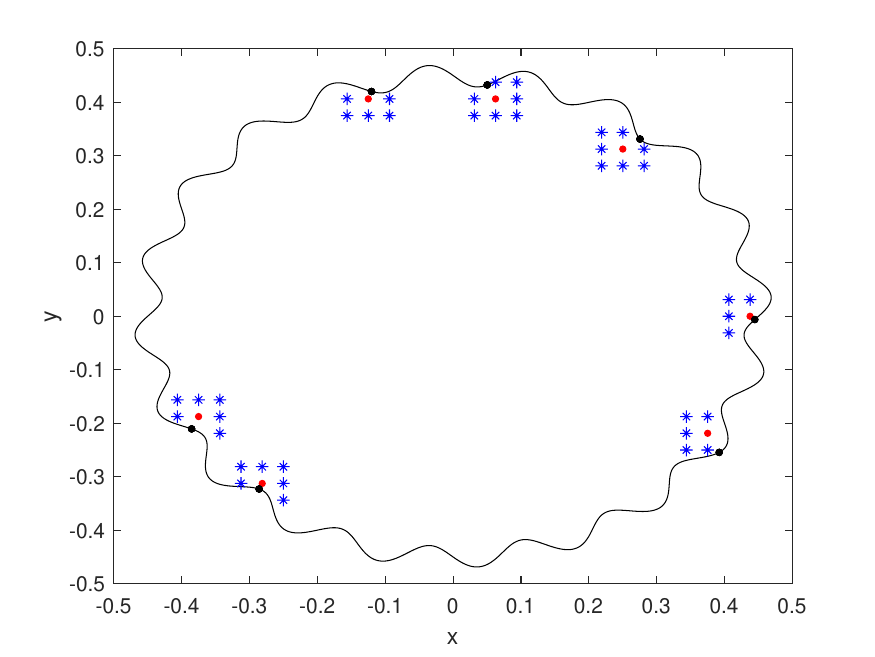}
	\end{subfigure}
	\caption
	{Illustrations for compact FDMs at irregular stencil centers, where each red point
is the stencil center of a boundary irregular stencil near $\partial \Omega$, a black point is the orthogonal projection of the stencil center of the boundary stencil onto $\partial \Omega$, the blue star points are remaining grid points used in compact FDMs, and the black curve is the boundary $\partial \Omega$.}
	\label{irregular:FDM:fig}
\end{figure}

	This paper is organized as following:
	In \cref{sec:FDM:regular}, we provide the fourth-order compact 9-point FDM at the regular  stencil center point.
	In \cref{sec:FDM:irregular}, we derive the fourth-order compact FDM at the irregular stencil center point in 10 cases.
	Since the  stencil of the compact FDM at the irregular  stencil center point has various configurations (see \cref{irregular:FDM:fig} for illustrations) for  the complex boundary curve, it is hard to derive the corresponding FDM for every configuration.
	 We present an at most 6 times 24 linear system with the explicit expression to compute the left-hand side of the stencil of the fourth-order compact FDM for any irregular  stencil center  points. Once the left-hand side of the FDM is determined, the right-hand side is calculated directly by its explicit formula.
	In \cref{sec:Numeri}, we test 3 examples with complicated curved boundaries to verify convergence rates of errors in $l_2$ and $l_{\infty}$ norms of the proposed FDM.
	In \cref{sec:contribu}, the main contribution is summarized.

	\section{Fourth-order  compact FDM at the regular  stencil center point }\label{sec:FDM:regular}
	In this section,  we use the fourth-order compact 9-point FDM that proposed in \cite[Theorem 2.1] {Feng3D2026} to compute $u_h$ to approximate $u$ at the regular  stencil center point $(x_i,y_j)$.
	\begin{theorem}\label{thm:FDM:2D}
		Let $\alpha>0,\vec{\beta},\kappa, \phi,u$ be smooth in $\overline{\Omega}$ in \eqref{model:problem:2D}, $\partial \Omega$ be a smooth boundary curve, and functions $a,b,d,f$ be defined in \eqref{notation:abdf}. Then the following compact 9-point FDM \cite[Theorem 2.1] {Feng3D2026} (see the left panel of \cref{nine:point:orthog:fig} for illustrations)
		\be\label{FDMs:2D:u}
		\mathcal{L}_h (u_h)_{i,j} :=\frac{1}{h^2}\sum_{r,\ell=-1}^1 C_{r,\ell}\Big|_{(x,y)=(x_i,y_j)} (u_h)_{i+r,j+\ell}=F_{i,j}\Big|_{(x,y)=(x_i,y_j)},
		\ee
		achieves the fourth-order consistency for \eqref{model:problem:2D} at the regular  stencil center point $(x_i,y_j)$, where
		\begin{align}\label{C:Left:2D}
			 &C_{-1,-1}:=\tfrac{1}{6}-\tfrac{a+b}{12}h, \notag \\
			 &C_{-1,0}:=\tfrac{2}{3}-\tfrac{a}{3}h+\tfrac{1}{12}[a^2+ab+d+2a_x+a_y+b_x]h^2 \notag \\
			&\qquad \quad -\tfrac{1}{24}[a(a_x+d)+ba_y+2d_x+\Delta a]h^3, \notag \\
			 &C_{-1,1}:=\tfrac{1}{6}-\tfrac{1}{12}[a-b]h-\tfrac{1}{12}[ab+a_y+b_x]h^2,  \notag \\
			 &C_{0,-1}:=\tfrac{2}{3}-\tfrac{b}{3}h+\tfrac{1}{12}[b^2+ab+a_y+b_x+2b_y+d]h^2 \notag\\
			&\qquad \quad-\tfrac{1}{24}[ab_x+b(b_y+d)+2d_y+\Delta b]h^3+\tfrac{bd_y}{12}h^4, \notag\\
			 &C_{0,0}:=\tfrac{-10}{3}-\tfrac{1}{6} [ a^2+ab+b^2-4d+2a_x+a_y+b_x+2b_y ]h^2 \notag \\
			&\qquad \quad+\tfrac{1}{12}[ad_x+\Delta d]h^4,\\
			 &C_{0,1}:=\tfrac{2}{3}+\tfrac{b}{3}h+\tfrac{1}{12}[ab+b^2+d+a_y+b_x+2b_y]h^2 \notag \\
			&\qquad \quad+\tfrac{1}{24}[ ab_x+b(b_y+d)+2d_y+\Delta b ]h^3, \notag \\
			 &C_{1,-1}:=\tfrac{1}{6}+\tfrac{1}{12}[a-b]h-\tfrac{1}{12}[ab+a_y+b_x]h^2-\tfrac{bd_y}{12}h^4, \notag \\
			&C_{1,0}:= \tfrac{2}{3}+\tfrac{a}{3}h+\tfrac{1}{12}[a^2+ab+d+2a_x+a_y+b_x]h^2 \notag \\
			&\qquad \quad+\tfrac{1}{24}[a(a_x+d)+ba_y+2d_x+\Delta a]h^3+\tfrac{bd_y}{12}h^4, \notag \\
			&C_{1,1}:=\tfrac{1}{6} +\tfrac{1}{12}[a+b]h, \notag
		\end{align}
		and
		\be\label{F:Right:2D}
		F_{i,j}:=f+\tfrac{1}{12}[ af_x+bf_y+\Delta f  ]h^2.
		\ee
	\end{theorem}
	By \eqref{notation:abdf}, \eqref{C:Left:2D}, and \eqref{F:Right:2D}, the first-order and second-order partial derivatives of $\alpha,\vec{\beta},\kappa, \phi$ in \eqref{model:problem:2D} are sufficient to attain the fourth-order consistency of the compact 9-point FDM at the regular  stencil center point in \cref{thm:FDM:2D}.
	
	\section{Fourth-order  compact FDM at an irregular  stencil center point }\label{sec:FDM:irregular}
	
	In this section,  we propose the fourth-order compact FDM at an irregular  stencil center point. For simplification, we assume that the parametric equation of the boundary curve $\partial\Omega$ satisfies
	\be\label{xtyt}
	(x,y)=(x(t),y(t))=(p(t),q(t)),\qquad \text{for }(x,y)\in \partial \Omega.
	\ee
	Then the Dirichlet boundary condition in \eqref{model:problem:2D} results in
	\be\label{Dirichlet}
	u(x(t),y(t))=g(t), \qquad (x(t),y(t))\in \partial \Omega.
	\ee
	For the irregular  stencil center point $(x_i,y_j)$, we define that
	\be\label{umn}
	 u^{(m,n)}:=\frac{\partial^{m+n}u(x,y)}{\partial x^m \partial y^n}\Big|_{(x,y)=({x}_i^o,{y}_j^o)},
	\ee
	\be\label{abmn}
	 a^{(m,n)}:=\frac{\partial^{m+n}a(x,y)}{\partial x^m \partial y^n}\Big|_{(x,y)=({x}_i^o,{y}_j^o)},\qquad b^{(m,n)}:=\frac{\partial^{m+n}b(x,y)}{\partial x^m \partial y^n}\Big|_{(x,y)=({x}_i^o,{y}_j^o)},
	\ee
	\be\label{dfmn}
	 d^{(m,n)}:=\frac{\partial^{m+n}d(x,y)}{\partial x^m \partial y^n}\Big|_{(x,y)=({x}_i^o,{y}_j^o)},\qquad f^{(m,n)}:=\frac{\partial^{m+n}f(x,y)}{\partial x^m \partial y^n}\Big|_{(x,y)=({x}_i^o,{y}_j^o)},
	\ee
	where $a,b,d,f$ are defined in \eqref{notation:abdf}, and $({x}_i^o,{y}_j^o)$ is the orthogonal projection on $\partial \Omega$ of the irregular  stencil center point $(x_i,y_j)$  (see black points in \cref{nine:point:orthog:fig} and \cref{irregular:FDM:fig} for illustrations).
	
	We also assume that there exists $ t_k^o$ such that
	\be\label{xiyj:curve}
	 ({x}_i^o,{y}_j^o)=(p(t_k^o),q(t_k^o)) \quad \text{for some}\ k \quad \text{with} \quad [p'(t_k^o)]^2+[q'(t_k^o)]^2\ne 0,
	\ee
	and
	\be\label{pqmn}
	p^{(m)}:=\frac{d^{m}p(t)}{d t^m }\Big|_{t=t_k^o},\qquad 	 q^{(m)}:=\frac{d^{m}q(t)}{d t^m }\Big|_{t=t_k^o}, \qquad g^{(m)}:=\frac{d^{m}g(t)}{d t^m }\Big|_{t=t_k^o}.
	\ee

		From our high-order FDM on the rectangle domain in \cite{FengMichelle2023,FengHanNeilan2026}, we need to use vertical and horizontal transformations for vertical  and horizontal boundaries, respectively.
	As $q^{(1)}/p^{(1)}$  is the slope of the tangent line of the boundary curve $\partial\Omega$, the tangent line is close to the vertical boundary case if $|q^{(1)}/p^{(1)}| \ge 1$, and the tangent line is close to the horizontal boundary case if $|q^{(1)}/p^{(1)}|\le 1$. So 
	we use the vertical transformation for $|q^{(1)}/p^{(1)}| \ge 1$ and the horizontal transformation for $|q^{(1)}/p^{(1)}|\le 1$ naturally.

	Based on the preliminary analysis,  to propose the fourth-order compact FDM at all irregular  stencil center points, we need to consider 10 subcases of 2 cases ($|q^{(1)}/p^{(1)}| \ge 1$  and $|q^{(1)}/p^{(1)}| \le 1$). In this section, we provide the derivation of FDM for the subcase 1 of the case 1 in \cref{thm:FDM:irregular:case:1}, and present the FDM for the subcase 1 of the case 2 directly in \cref{thm:FDM:irregular:case:2}. For other 8 subcases, we give necessary details to help readers to derive corresponding FDMs straightforwardly. Furthermore, we also present \cref{remark:4:subcases} to provide a feasible way to avoid 8 subcases, i.e., the FDMs in \cref{thm:FDM:irregular:case:1,thm:FDM:irregular:case:2} can achieve the fourth-order consistency at any irregular  stencil center points of all 10 subcases.
	
	\textbf{Case 1: $|q^{(1)}/p^{(1)}| \ge 1$.}
	\eqref{2D:linear:simple:eq} yields
	\be\label{f00}
	au^{(1,0)} +bu^{(0,1)} +du+	 u^{(2,0)} + u^{(0,2)}=f.
	\ee
	Differentiating \eqref{f00} with respect to $x$ and $y$, we obtain
	\be\label{f10}
	\begin{split}
		&au^{(2, 0)}+a^{(1, 0)}u^{(1, 0)}+bu^{(1, 1)}+b^{(1, 0)}u^{(0, 1)}+du^{(1, 0)}+d^{(1, 0)}u+u^{(1, 2)}+u^{(3, 0)}=f^{(1, 0)},
	\end{split}
	\ee
	\be\label{f01}
	\begin{split}
		&au^{(1, 1)}+a^{(0, 1)}u^{(1, 0)}+bu^{(0, 2)}+b^{(0, 1)}u^{(0, 1)}+du^{(0, 1)}+d^{(0, 1)}u+u^{(0, 3)}+u^{(2, 1)}=f^{(0, 1)}.
	\end{split}
	\ee
	Then \eqref{f00}--\eqref{f01} imply (see the first arrow in \cref{fig:case:1} for an illustration)
	\begin{align}\label{u20:u30}
		u^{(2, 0)} =& -au^{(1, 0)}-bu^{(0, 1)}-du-u^{(0, 2)}+f, \notag \\
		u^{(2, 1)} = & -au^{(1, 1)}-a^{(0, 1)}u^{(1, 0)}-bu^{(0, 2)}-b^{(0, 1)}u^{(0, 1)}-d u^{(0, 1)} \notag \\
		&-d^{(0, 1)}u-u^{(0, 3)}+f^{(0, 1)}, \\
		u^{(3, 0)} = & a^2u^{(1, 0)}+abu^{(0, 1)}+adu-af+au^{(0, 2)}-a^{(1, 0)}u^{(1, 0)}  \notag \\
		&-bu^{(1, 1)}-b^{(1, 0)}u^{(0, 1)}-du^{(1, 0)}-d^{(1, 0)}u-u^{(1, 2)}+f^{(1, 0)}. \notag
	\end{align}
	Differentiating \eqref{Dirichlet} with respect to $t$ three times yields
	\begin{align} \label{g0:g3}
		& u=g, \notag \\
		& p^{(1)}u^{(1, 0)}+q^{(1)}u^{(0, 1)}=g^{(1)}, \notag \\
		& [p^{(1)}]^2u^{(2, 0)}+2p^{(1)}q^{(1)}u^{(1, 1)}+[q^{(1)}]^2u^{(0, 2)} \notag  \\
		&\quad +p^{(2)}u^{(1, 0)}+q^{(2)}u^{(0, 1)}=g^{(2)},\notag \\
		& [p^{(1)}]^3u^{(3, 0)}+3[p^{(1)}]^2q^{(1)}u^{(2, 1)}+3p^{(1)}[q^{(1)}]^2u^{(1, 2)}\\
		&\quad +[q^{(1)}]^3u^{(0, 3)}+3p^{(1)}p^{(2)}u^{(2, 0)}  +3p^{(1)}q^{(2)}u^{(1, 1)}\notag \\
		&\quad +3p^{(2)}q^{(1)}u^{(1, 1)}+3q^{(1)}q^{(2)}u^{(0, 2)}+p^{(3)}u^{(1, 0)}\notag  \\
		&\quad +q^{(3)}u^{(0, 1)}=g^{(3)}. \notag
	\end{align}
	Substituting \eqref{u20:u30} into \eqref{g0:g3}, then the symbolic calculation results in (see the second arrow in \cref{fig:case:1} for an illustration)
	\begin{align} \label{u10:u12}
		u =& g, \notag \\
		u^{(0, 1)} =& s_{1}u^{(1, 0)}-z_{1}g^{(1)},\notag \\
		u^{(0, 2)} =& s_{2}u^{(1, 0)}+s_{3}u^{(1, 1)}-z_{3}f-z_{2}g-z_{4}g^{(1)}-z_{5}g^{(2)},\\
		u^{(0, 3)} =& s_{4}u^{(1, 0)}+s_{5}u^{(1, 1)}+s_{6}u^{(1, 2)}	 -z_{7}f-z_{8}f^{(0, 1)}\notag \\
		&-z_{10}f^{(1, 0)}-z_{6}g-z_{9}g^{(1)}-z_{11}g^{(2)}-z_{12}g^{(3)}, \notag
	\end{align}
	where
	\begin{align}\label{s1:s6}
		s_{1} =& -p^{(1)}/ q^{(1)}, \notag \\
		s_{2} =& \big\{ b[p^{(1)}]^3-a[p^{(1)}]^2q^{(1)}-p^{(1)}q^{(2)}\notag \\ & +p^{(2)}q^{(1)}\big\}/ \big \{q^{(1)}([p^{(1)}]^2-[q^{(1)}]^2)\big \}, \notag \\
		s_{3} = & 2q^{(1)}p^{(1)}/\big \{[p^{(1)}]^2-[q^{(1)}]^2 \big \}, \notag \\
		s_{4} =& \big\{ b^{(1, 0)}[p^{(1)}]^6-q^{(1)}(3b^2+a^{(1, 0)}-3b^{(0, 1)}\notag \\
		&-2d)[p^{(1)}]^5+([4ab-3a^{(0, 1)} -b^{(1, 0)}][q^{(1)}]^2\notag \\
		&-aq^{(2)})[p^{(1)}]^4+([a^{(1, 0)}-a^2-3b^{(0, 1)}-2d][q^{(1)}]^3\notag \\
		 &+(ap^{(2)}+6bq^{(2)})q^{(1)}-q^{(3)})[p^{(1)}]^3+(3a^{(0, 1)}[q^{(1)}]^4\notag \\
		 &-(3aq^{(2)}+6bp^{(2)})[q^{(1)}]^2+p^{(3)}q^{(1)}+3p^{(2)}q^{(2)})[p^{(1)}]^2\\
		 &+3(ap^{(2)}[q^{(1)}]^2-[p^{(2)}]^2+q^{(1)}q^{(3)}/3-[q^{(2)}]^2)q^{(1)}p^{(1)}\notag \\
		 &+3p^{(2)}q^{(2)}[q^{(1)}]^2-p^{(3)}[q^{(1)}]^3\big\}/\big \{ [q^{(1)}]^2(3[p^{(1)}]^4\notag \\
		 &-4[p^{(1)}]^2[q^{(1)}]^2+[q^{(1)}]^4) \big\}, \notag \\
		s_{5} = &\big\{ 3a[p^{(1)}]^2[q^{(1)}]^3-a[p^{(1)}]^4q^{(1)}-b[p^{(1)}]^5-5b[p^{(1)}]^3[q^{(1)}]^2\notag \\
		 &+3[p^{(1)}]^3q^{(2)}-3[p^{(1)}]^2p^{(2)}q^{(1)}+3p^{(1)}[q^{(1)}]^2q^{(2)}-3p^{(2)}[q^{(1)}]^3\big\}\notag \\
		&/\big \{q^{(1)}(3[p^{(1)}]^4-4[p^{(1)}]^2[q^{(1)}]^2+[q^{(1)}]^4) \big \}, \notag \\
		s_{6} =& -p^{(1)}\big\{[p^{(1)}]^2-3[q^{(1)}]^2\big\}/\big \{q^{(1)}(3[p^{(1)}]^2-[q^{(1)}]^2) \big\}, \notag
	\end{align}
	and
	\begin{align}\label{z1:z12}
		z_{1} =& -1/q^{(1)}, \notag \\
		z_{2} =& d[p^{(1)}]^2/\big\{ [p^{(1)}]^2-[q^{(1)}]^2 \big\}, \notag \\
		z_{3} =& -[p^{(1)}]^2/\big\{ [p^{(1)}]^2-[q^{(1)}]^2 \big\}, \notag \\
		z_{4} =& \big\{ b[p^{(1)}]^2-q^{(2)}\big\}/\big\{ q^{(1)}([p^{(1)}]^2-[q^{(1)}]^2)\big\}, \notag \\
		z_{5} =& 1/\big\{[p^{(1)}]^2-[q^{(1)}]^2\big\}, \notag \\
		z_{6} =& p^{(1)}\big\{ d^{(1, 0)}[p^{(1)}]^4-3q^{(1)}(bd-d^{(0, 1)})[p^{(1)}]^3+[q^{(1)}]^2(ad\notag \\
		&-d^{(1, 0)})[p^{(1)}]^2 +(3dq^{(1)}q^{(2)}-3d^{(0, 1)}[q^{(1)}]^3)p^{(1)}-3dp^{(2)}[q^{(1)}]^2\big\} \notag \\
		 &/\big\{q^{(1)}(3[p^{(1)}]^4-4[p^{(1)}]^2[q^{(1)}]^2+[q^{(1)}]^4)\big\}, \notag \\
		z_{7} = &-p^{(1)}\big\{ a[p^{(1)}]^2q^{(1)}-3b[p^{(1)}]^3+3p^{(1)}q^{(2)}-3p^{(2)}q^{(1)}\big\}/\big\{ 3[p^{(1)}]^4\notag \\
		&-4[p^{(1)}]^2[q^{(1)}]^2 +[q^{(1)}]^4\big\}, \notag \\
		z_{8} =& -3[p^{(1)}]^2/\big\{ 3[p^{(1)}]^2-[q^{(1)}]^2\big\}, \\
		z_{9} =& \big\{ b^{(1, 0)}[p^{(1)}]^5-3q^{(1)}(b^2-b^{(0, 1)}-d)[p^{(1)}]^4+([ab-b^{(1, 0)}][q^{(1)}]^2 \notag \\
		& -aq^{(2)})[p^{(1)}]^3+(6bq^{(1)}q^{(2)}-3(d+b^{(0, 1)})[q^{(1)}]^3-q^{(3)})[p^{(1)}]^2\notag \\
		 &-3p^{(2)}(b[q^{(1)}]^2-q^{(2)})p^{(1)}+q^{(3)}[q^{(1)}]^2-3q^{(1)}[q^{(2)}]^2\big\}\notag \\
		&/\big\{ [q^{(1)}]^2 (3[p^{(1)}]^4-4[p^{(1)}]^2[q^{(1)}]^2+[q^{(1)}]^4)\big\}, \notag\\
		z_{10} =& -[p^{(1)}]^3/\big\{ q^{(1)}(3[p^{(1)}]^2-[q^{(1)}]^2)\big\}, \notag \\
		z_{11} =& \big\{a[p^{(1)}]^3-3b[p^{(1)}]^2q^{(1)}-3p^{(1)}p^{(2)}+3q^{(1)}q^{(2)}\big\}/\big\{ q^{(1)}(3[p^{(1)}]^4\notag \\
		 &-4[p^{(1)}]^2[q^{(1)}]^2+[q^{(1)}]^4)\big\}, \notag \\
		z_{12} =& 1/\big\{ q^{(1)}(3[p^{(1)}]^2-[q^{(1)}]^2) \big\}. \notag
	\end{align}
	From denominators of expressions in \eqref{s1:s6} and \eqref{z1:z12}, we can say that  $s_{1},\dots s_{6}$ and $z_1\dots z_{12}$  in \eqref{s1:s6}  and \eqref{z1:z12} are well defined if
	\be\label{well:defined:1}
	q^{(1)}\ne 0,\qquad q^{(1)}\ne \pm p^{(1)},\qquad  q^{(1)}\ne \pm\sqrt{3} p^{(1)}.
	\ee
	As we are considering the case 1 with $|q^{(1)}/p^{(1)}| \ge 1$, $q^{(1)}=0$ can be ignored for the case 1, which is also the reason that we choose the vertical transformation in the first transformation  (i.e., the first arrow) in \cref{fig:case:1} for  $|q^{(1)}/p^{(1)}| \ge 1$. So we discuss following 5 subcases of the case 1:

	\textbf{Subcase 1 of the case 1: $q^{(1)}\neq \pm p^{(1)}$ and   $q^{(1)}\neq \pm\sqrt{3} p^{(1)}$.}
	\begin{figure}[h]
		 \begin{subfigure}[b]{0.3\textwidth}
			\hspace{2cm}
			\begin{tikzpicture}[scale = 0.45]
				\node (A) at (0,7) {$u^{(0,0)}$};
				\node (A) at (0,6) {$u^{(0,1)}$};
				\node (A) at (0,5) {$u^{(0,2)}$};
				\node (A) at (0,4) {$u^{(0,3)}$};
				\node (A) at (2,7) {$u^{(1,0)}$};
				\node (A) at (2,6) {$u^{(1,1)}$};
				\node (A) at (2,5) {$u^{(1,2)}$};
				\node (A) at (4,7) {$u^{(2,0)}$};
				\node (A) at (4,6) {$u^{(2,1)}$};
				\node (A) at (6,7) {$u^{(3,0)}$};
				\draw[->, >=stealth, line width=1pt, scale=1.5] (3.2,3.5) -- (7,3.5);
			\end{tikzpicture}
		\end{subfigure}
		 \begin{subfigure}[b]{0.3\textwidth}
			\hspace{2cm}
			\begin{tikzpicture}[scale = 0.45]
				\node (A) at (0,7) {$u^{(0,0)}$};
				\node (A) at (0,6) {$u^{(0,1)}$};
				\node (A) at (0,5) {$u^{(0,2)}$};
				\node (A) at (0,4) {$u^{(0,3)}$};
				\node (A) at (2,7) {$u^{(1,0)}$};
				\node (A) at (2,6) {$u^{(1,1)}$};
				\node (A) at (2,5) {$u^{(1,2)}$};
				\draw[->, >=stealth, line width=1pt, scale=1.5] (2.2,3.5) -- (6,3.5);	  	
			\end{tikzpicture}
		\end{subfigure}
		 \begin{subfigure}[b]{0.3\textwidth}
			\hspace{1.2cm}
			\vspace{0.4cm}
			\begin{tikzpicture}[scale = 0.45]
				\node (A) at (0,7) {$u^{(1,0)}$};
				\node (A) at (0,6) {$u^{(1,1)}$};
				\node (A) at (0,5) {$u^{(1,2)}$};
			\end{tikzpicture}
		\end{subfigure}	
		\caption
		{The illustration of vertical transformations for the subcase 1 with $q^{(1)}\neq \pm p^{(1)}$ and   $q^{(1)}\neq \pm\sqrt{3} p^{(1)}$ of the case 1 with $|q^{(1)}/p^{(1)}| \ge 1$. }
		\label{fig:case:1}
	\end{figure}		
	Recall that the Taylor series generates that	
	\be \label{Taylor:Se}
	 u(x+{x}_i^o,y+{y}_j^o)=\sum_{\substack{m,n=0 \\ m+n\le 3}}^3
	\frac{x^m y^n}{m! n!}	 u^{(m,n)}+\bo(h^{4})\quad \text{with} \quad (x,y)\in (-2h,2h)^2,
	\ee
	where $u^{(m,n)}$ is defined in \eqref{umn} and \eqref{xiyj:curve}.  Plugging \eqref{u20:u30} and \eqref{u10:u12} into \eqref{Taylor:Se}, we have
	\be \label{Taylor:New}
	\begin{split}
		u(x+{x}_i^o,y+{y}_j^o)=& \sum_{n=0}^2 G_{1,n}(x,y)u^{(1,n)}+\sum_{\substack{m,n=0 \\ m+n\le 1}}^1 H_{m,n}(x,y) f^{(m,n)}\\
		&+ \sum_{n=0}^3V_n(x,y)g^{(n)}+\bo(h^{4}) \quad \text{with} \quad (x,y)\in (-2h,2h)^2,
	\end{split}
	\ee
	where
	\begin{align} \label{GH}
		& G_{1,0}(x,y)=x+ \xi_{1}y+ \xi_{2}y^2+ \xi_{3}y^3+\xi_{4}x^2 +\xi_{5}x^2y +\xi_{6}x^3, \notag \\
		& G_{1,1}(x,y)=\xi_{7}y^2 + \xi_{8}y^3+ \xi_{9}x^2+\xi_{10}x^2y + \xi_{11}x^3+xy, \notag\\
		& G_{1,2}(x,y)=\xi_{12}y^3+ \xi_{13}x^2y-x^3/6+xy^2/2, \notag \\
		& H_{0,0}(x,y)= \eta_{1} y^2 +\eta_{2} y^3 +\eta_{3} x^2 +\eta_{4} x^2y + \eta_{5} x^3, \notag \\
		&H_{1,0}(x,y)= \eta_{8}y^3+ \eta_{9}x^2y+x^3/6, \quad H_{0,1}(x,y)=\eta_{6} y^3 + \eta_{7} x^2y,   \\
		& V_{0}(x,y)=  1+\omega_{1} y^2 +\omega_{2} y^3 +\omega_{3} x^2 + \omega_{4} x^2y+\omega_{5}x^3,\notag \\
		& V_{1}(x,y)= \omega _{6}y +\omega _{7} y^2  +\omega _{8} y^3 +\omega _{9} x^2 +\omega _{10} x^2y + \omega _{11} x^3,\notag \\
		& V_{2}(x,y)= \omega_{12} y^2+\omega_{13} y^3 + \omega_{14} x^2+\omega_{15} x^2y +\omega_{16} x^3 ,\notag \\
		& V_{3}(x,y)= \omega_{17}y^3 + \omega_{18}x^2y, \notag
	\end{align}	
	with
	\begin{align}\label{xi:eta:ome}
	&	\xi_{1} = s_{1}, \qquad  \xi_{2} = s_{2}/2, \qquad  \xi_{3} = s_{4}/6, \qquad  \xi_{4} = -(bs_{1}+a+s_{2})/2, \notag \\
	&  \xi_{5} = -[(d+b^{(0, 1)})s_{1}+s_{2}b+s_{4}+a^{(0, 1)}]/2, \notag \\
	&  \xi_{6} = (a^2+(bs_{1}+s_{2})a-s_{1}b^{(1, 0)}-a^{(1, 0)}-d)/6, \qquad  \xi_{7} = s_{3}/2, \notag \\
	&  \xi_{8} = s_{5}/6, \qquad  \xi_{9} = -s_{3}/2, \qquad  \xi_{10} = -(s_{3}b+a+s_{5})/2,   \notag \\
	& \xi_{11} = (s_{3}a-b)/6, \qquad  \xi_{12} = s_{6}/6, \qquad  \xi_{13} = -s_{6}/2, \notag \\
	& \eta_{1} = -z_{3}/2, \qquad  \eta_{2} = -z_{7}/6, \qquad  \eta_{3} = (z_{3}+1)/2, \notag  \\
	&  \eta_{4} = (z_{3}b+z_{7})/2, \qquad \eta_{5} = -(z_{3}+1)a/6, \qquad  \eta_{6} = -z_{8}/6,    \notag \\
	& \eta_{7} = (z_{8}+1)/2, \qquad \eta_{8} = -z_{10}/6, \qquad  \eta_{9} = z_{10}/2, \\
	&	\omega_{1} = -z_{2}/2, \qquad  \omega_{2} = -z_{6}/6, \qquad  \omega_{3} = (z_{2}-d)/2,  \notag \\ & \omega_{4} = (z_{2}b-d^{(0, 1)}+z_{6})/2, \qquad \omega_{5} = [a(d-z_{2})-d^{(1, 0)}]/6,  \notag \\
	&   \omega_{6} = -z_{1}, \qquad  \omega_{7} = -z_{4}/2, \qquad  \omega_{8} = -z_{9}/6, \notag \\
	&  \omega_{9} = (bz_{1}+z_{4})/2, \qquad  \omega_{10} = [z_{1}(b^{(0, 1)}+d)+z_{4}b+z_{9}]/2,  \notag \\
	&  \omega_{11} = [b^{(1, 0)}z_{1}-a(bz_{1}+z_{4})]/6, \qquad  \omega_{12} = -z_{5}/2,     \notag \\
	&   \omega_{13} = -z_{11}/6, \qquad  \omega_{14} = z_{5}/2, \qquad \omega_{15} = (z_{5}b+z_{11})/2,   \notag \\
	&   \omega_{16} = -z_{5}a/6, \qquad  \omega_{17} = -z_{12}/6, \qquad  \omega_{18} = z_{12}/2,\notag
	\end{align}
	and $s_{1},\dots s_{6}$ and $z_1\dots z_{12}$ are defined in \eqref{s1:s6}  and \eqref{z1:z12}, respectively.
	Let define
	\be\label{Crl}
	C_{r,\ell}:=\sum_{\sigma =0}^3 c_{r,\ell,\sigma} h^\sigma, \quad \text{with }  c_{r,\ell,\sigma} \text{ depending on } a,b,d \text{ in } \eqref{2D:linear:simple:eq}, \quad \text{and} \quad c_{r,\ell,\sigma}|_{(x,y)=({x}_i^o,{y}_j^o)} \in \R,
	\ee
	and assume that
	\be \label{set:S}
	\begin{split}
		& S:=\{ (r,\ell) :   ({x}_i^o+rh,{y}_j^o+\ell h)\in \overline{\Omega} \text{ with } r,\ell=-1,0,1  \},\\
		& w=(x_i-{x}_i^o)/h,\qquad  v=(y_j-{y}_j^o)/h,
	\end{split}
	\ee
	where $({x}_i^o,{y}_j^o)\in \Omega$ is the orthogonal projection of the irregular  stencil center point $(x_i,y_j)$ (the set $S$ consists of blue star and red points in \cref{irregular:FDM:fig}). Next, we propose the derivation of the fourth-order compact FDM at the irregular  stencil center point.

	{\bf{ The algorithm to derive the fourth-order compact FDM at an irregular  stencil center:}}	We define
	\be\label{Lh:u:1}
	\begin{split}
		\mathcal{L}_h u_{i,j} :=&\sum_{\substack{r,\ell=-1 \\ (r,\ell)\in S}}^1 C_{r,\ell}\Big|_{(x,y)=({x}_i^o,{y}_j^o)} u_{i+r,j+\ell}\\
		=&\sum_{\substack{r,\ell=-1 \\ (r,\ell)\in S}}^1 C_{r,\ell}\Big|_{(x,y)=({x}_i^o,{y}_j^o)} u(rh+wh+{x}_i^o,\ell h+vh+{y}_j^o),
	\end{split}
	\ee
	where $S,w,v$ are defined in \eqref{set:S}.
	Applying \eqref{Taylor:New} into \eqref{Lh:u:1}, we obtain
	\begin{align}\label{Lh:u}
		\mathcal{L}_h u_{i,j}   = &	  \sum_{\substack{r,\ell=-1 \\ (r,\ell)\in S}}^1 C_{r,\ell} \Big|_{({x}_i^o,{y}_j^o)}	 \sum_{n=0}^2
		G_{1,n} ([r+w]h,[\ell+v]h) u^{(1,n)}\notag \\
		&+ \sum_{\substack{r,\ell=-1 \\ (r,\ell)\in S}}^1 C_{r,\ell} \Big|_{({x}_i^o,{y}_j^o)} \sum_{\substack{m,n=0 \\ m+n\le 1}}^1 H_{m,n} ([r+w]h,[\ell+v]h)	f^{(m,n)}  \\
		& +\sum_{\substack{r,\ell=-1 \\ (r,\ell)\in S}}^1 C_{r,\ell} \Big|_{({x}_i^o,{y}_j^o)}\sum_{n=0}^3V_n ([r+w]h,[\ell+v]h) g^{(n)} +\bo(h^{4}).\notag
	\end{align}
	Let
	\be\label{Lh:uh:1}
	\begin{split}
		\mathcal{L}_h (u_h)_{i,j} :=&\sum_{\substack{r,\ell=-1 \\ (r,\ell)\in S}}^1 C_{r,\ell}\Big|_{(x,y)=({x}_i^o,{y}_j^o)} (u_h)_{i,j}=F_{i,j},
	\end{split}
	\ee
	where
	\begin{align}\label{Fij}
		F_{i,j}:=\text{the terms of } \Bigg( & \sum_{\substack{r,\ell=-1 \\ (r,\ell)\in S}}^1 C_{r,\ell} \Big|_{({x}_i^o,{y}_j^o)} \sum_{\substack{m,n=0 \\ m+n\le 1}}^1 H_{m,n} ([r+w]h,[\ell+v]h) 	f^{(m,n)} \notag \\
		& +\sum_{\substack{r,\ell=-1 \\ (r,\ell)\in S}}^1 C_{r,\ell} \Big|_{({x}_i^o,{y}_j^o)}\sum_{n=0}^3V_n ([r+w]h,[\ell+v]h) g^{(n)} \Bigg ) \\
		& \text{with degree}   \le 3 \text{ in } h. \notag
	\end{align}
	Then
	\be\label{L:h:u:uh}
	\mathcal{L}_h (u-u_h)_{i,j}=	 \sum_{\substack{r,\ell=-1 \\ (r,\ell)\in S}}^1 C_{r,\ell} \Big|_{({x}_i^o,{y}_j^o)}	 \sum_{n=0}^2
	G_{1,n} ([r+w]h,[\ell+v]h)  u^{(1,n)} +\bo(h^{4}).
	\ee
	To achieve the fourth-order consistency, we need to solve
	\be\label{oh4:final}
	\sum_{\substack{r,\ell=-1 \\ (r,\ell)\in S}}^1 C_{r,\ell} \Big|_{({x}_i^o,{y}_j^o)}	 \sum_{n=0}^2
	G_{1,n} ([r+w]h,[\ell+v]h)  u^{(1,n)} =\bo(h^{4}).
	\ee
	By definitions of $C_{k,\ell}$ in \eqref{Crl} and $	G_{1,n}(x,y)$ in \eqref{GH}, we have that \eqref{oh4:final} is equivalent to solve the following linear system
	
	\be \label{sub:linear:system}
	\sum_{\substack{r,\ell=-1 \\ (r,\ell)\in S}}^1 	A_{r,\ell} c_{r,\ell}=0,\quad
	c_{r,\ell}:=\begin{bmatrix}
		c_{r,\ell,0}\\
		c_{r,\ell,1}\\
		c_{r,\ell,2}\\
		c_{r,\ell,3}
	\end{bmatrix},
	\ee
	with the $6 \times 4$ lower triangular matrix $A_{r,\ell}$
		\be\label{maxtrix:1}
	A_{r,\ell}:=\begin{bmatrix}
\tau_{\ell}\xi_{1}+\mu_{r} & 0 & 0 &  0\\
\mu_{r}^2\xi_{9}+\tau_{\ell}^2\xi_{7}+\mu_{r}\tau_{\ell} & 0 & 0 & 0\\
\mu_{r}^2\tau_{\ell}\xi_{13}+\tau_{\ell}^3\xi_{12}-\mu_{r}^3/6+\mu_{r}\tau_{\ell}^2/2 & 0 & 0 & 0\\
\mu_{r}^2\xi_{4}+\tau_{\ell}^2\xi_{2} & \tau_{\ell}\xi_{1}+\mu_{r} & 0 &  0\\
\mu_{r}^3\xi_{11}+\mu_{r}^2\tau_{\ell}\xi_{10}+\tau_{\ell}^3\xi_{8} & \mu_{r}^2\xi_{9}+\tau_{\ell}^2\xi_{7}+\mu_{r}\tau_{\ell} & 0 & 0\\
\mu_{r}^3\xi_{6}+\mu_{r}^2\tau_{\ell}\xi_{5}+\tau_{\ell}^3\xi_{3} &  \mu_{r}^2\xi_{4}+\tau_{\ell}^2\xi_{2} & \tau_{\ell}\xi_{1}+\mu_{r} & 0
	\end{bmatrix},
	\ee
	where
	\be\label{mur:tau:l}
	\mu_{r}=r+w, \qquad \tau_{\ell}=\ell+v,
	\ee
$\xi_{1},\dots \xi_{13}$ are defined in  \eqref{xi:eta:ome},  $S,w$ and $v$ are defined in \eqref{set:S}, $s_{1},\dots s_{6}$ and $z_1\dots z_{12}$ are defined in \eqref{s1:s6}  and \eqref{z1:z12}. From \eqref{maxtrix:1}, the fourth column of $A_{r,\ell}$ is the zero vector, so we can set all $c_{r,\ell,3}=0$ in \eqref{sub:linear:system}.
As we have at most 8 grid points in the stencil of the compact FDM at the irregular  stencil center point, the matrix size of the linear system \eqref{sub:linear:system} is at most $6$ by $(4-1)8$ = $6$ by $24$.
	
	Substituting $C_{k,\ell}$ of \eqref{Crl}, $H_{m,n}(x,y)$ and $V_{n}(x,y)$ of \eqref{GH} into \eqref{Fij}, we have that the explicit expression of $F_{i,j}$ in \eqref{Fij} is
	\begin{align}  \label{Fij:explicit}
		 F_{i,j}=&\sum_{\substack{r,\ell=-1 \\ (r,\ell)\in S}}^1 \bigg\{c_{r, \ell, 0}g+ \Big \{c_{r, \ell, 0}g^{(1)}\omega_{6}\tau_{\ell}+c_{r, \ell, 1}g \Big \} h+\Big \{ c_{r, \ell, 0}\big [(\eta_{1}f+g\omega_{1}+g^{(1)}\omega_{7} \notag \\
		 &+g^{(2)}\omega_{12})\tau_{\ell}^2 +\mu_{r}^2(\eta_{3}f+g\omega_{3}+g^{(1)}\omega_{9}+g^{(2)}\omega_{14})\big ]+c_{r, \ell, 1}g^{(1)}\omega_{6}\tau_{\ell}+c_{r, \ell, 2}g\Big \}h^2 \notag \\
		&+\Big \{  c_{r, \ell, 0}\big[ (\eta_{2}f+\eta_{6}f^{(0, 1)}+\eta_{8}f^{(1, 0)}+g\omega_{2}+g^{(1)}\omega_{8}+g^{(2)}\omega_{13}+g^{(3)}\omega_{17})\tau_{\ell}^3 \notag \\
		 &+\mu_{r}^2(\eta_{4}f+\eta_{7}f^{(0, 1)}+\eta_{9}f^{(1, 0)}+g\omega_{4}+g^{(1)}\omega_{10}+g^{(2)}\omega_{15}+g^{(3)}\omega_{18})\tau_{\ell}  \\
		 &+\mu_{r}^3(\eta_{5}f+\omega_{5}g+\omega_{11}g^{(1)}+\omega_{16}g^{(2)}+f^{(1, 0)}/6)\big]+c_{r, \ell, 1}(\eta_{1}f +g\omega_{1}+g^{(1)}\omega_{7}\notag \\
		 &+g^{(2)}\omega_{12})\tau_{\ell}^2+c_{r, \ell, 2}g^{(1)}\omega_{6}\tau_{\ell}+c_{r, \ell, 1}(\eta_{3}f+g\omega_{3}+g^{(1)}\omega_{9} +g^{(2)}\omega_{14})\mu_{r}^2+c_{r, \ell, 3}g \Big \} h^3\bigg\}, \notag
	\end{align}
	$\mu_{r}$ and $\tau_{\ell}$ are defined in \eqref{mur:tau:l},	 $\eta_{1},\dots \eta_{9}$, and  $\omega_{1},\dots \omega_{18}$ are defined in  \eqref{xi:eta:ome},  $S,w$ and $v$ are defined in \eqref{set:S}, $s_{1},\dots s_{6}$ and $z_1\dots z_{12}$ are defined in \eqref{s1:s6}  and \eqref{z1:z12}, $g^{(1)},g^{(2)}$, and $g^{(3)}$ are defined in \eqref{pqmn}.
	
	Now, we can say that the compact FDM in \eqref{Lh:uh:1} approximates the Dirichlet boundary condition	$u =g$ at the point $(x_i^o,y_j^o)\in \partial \Omega$ with the fourth-order consistency, if 	 $C_{r,\ell}=\sum_{\sigma=0}^3 c_{r,\ell,\sigma} h^\sigma$  and $F_{i,j}$ in \eqref{Fij:explicit} are obtained by  solving 	 \eqref{sub:linear:system} with \eqref{maxtrix:1} and satisfying $c_{0,0,0}\ne 0$ and $F_{i,j}|_{h=0}=g$. Using \eqref{Fij:explicit}, $F_{i,j}|_{h=0}=g$ if and only if $\sum_{\substack{r,\ell=-1 \\ (r,\ell)\in S}}^1 c_{r, \ell, 0}=1$.

	In summary, we propose the fourth-order compact FDM to solve \eqref{model:problem:2D} at the irregular  stencil center point for the subcase 1 of the case 1 as follows:
	 \begin{theorem}\label{thm:FDM:irregular:case:1}
		Let $\alpha>0,\vec{\beta},\kappa, \phi,u$ be smooth in $\overline{\Omega}$ in \eqref{model:problem:2D}, $\partial \Omega$ be a smooth boundary curve,  functions $a,b,d,f$ be defined in \eqref{notation:abdf}, and $({x}_i^o,{y}_j^o)\in \Omega$ be the orthogonal projection of the irregular  stencil center point $(x_i,y_j)$. If
		\[
			 ({x}_i^o,{y}_j^o)=(p(t_k^o),q(t_k^o)), \quad |q'(t_k^o)/p'(t_k^o)| \ge 1, \quad q'(t_k^o)\ne \pm p'(t_k^o), \quad q'(t_k^o)\ne \pm \sqrt{3} p'(t_k^o).
		\]
		Then the following compact FDM  (see  \cref{irregular:FDM:fig} for illustrations)
		\[
		\mathcal{L}_h (u_h)_{i,j} :=\sum_{\substack{r,\ell=-1 \\ (r,\ell)\in S}}^1 C_{r,\ell}\Big|_{(x,y)=(x_i^o,y_j^o)} (u_h)_{i+r,j+\ell}=F_{i,j}\Big|_{(x,y)=(x^o_i,y^o_j)},
		\]
		achieves the fourth-order consistency for the Dirichlet boundary condition	$u =g$ at the point $(x_i^o,y_j^o)\in \partial\Omega$, where $S$ is defined in \eqref{set:S}, the right-hand side $F_{i,j}$ of the stencil is defined in \eqref{Fij:explicit}, the left-hand side  $C_{r,\ell}=\sum_{\sigma=0}^3 c_{r,\ell,\sigma} h^\sigma$  of the stencil  is obtained by  solving 	 \eqref{sub:linear:system} with \eqref{maxtrix:1}, $c_{0,0,0}\ne 0$, and $\sum_{\substack{r,\ell=-1 \\ (r,\ell)\in S}}^1 c_{r, \ell, 0}=1$.
	\end{theorem}
	\begin{remark}
	From \eqref{u20:u30} and \eqref{u10:u12}--\eqref{z1:z12}, we observe that first-order derivatives of  $\alpha,\vec{\beta},\kappa, \phi$ in \eqref{model:problem:2D} and up to third-order derivatives of the Dirichlet boundary function $g$ and the parametric equation $(p(t),q(t))$ in \eqref{xtyt} are sufficient to construct the fourth-order compact FDM at the irregular  stencil center point.		
	\end{remark}
	\begin{remark}
		Note that the solution $c_{r,\ell}$ in	 \eqref{sub:linear:system} is not unique. In our numerical examples, we set $c_{0,0,0}=1$ and use the MATLAB function 'mldivide' to solve \eqref{sub:linear:system} to set the remaining free parameters to zero automatically.
	\end{remark}
	 \begin{remark}\label{remark:4:subcases}
	We still need to derive FDMs for the remaining 4 subcases for the case 1. Readers can choose other  $(\tilde{x}^o_i,\tilde{y}^o_j)\in \partial \Omega$ that
	 is near  $({x}_i^o,{y}_j^o)$ to replace $({x}_i^o,{y}_j^o)$   to avoid the remaining 4 subcases. Then readers can utilize the FDM of the subcase 1 to completely handle any irregular  stencil center points in the case 1. In our numerical examples, each of $({x}_i^o,{y}_j^o)$ is the orthogonal projection of the irregular  stencil center point $({x}_i,{y}_j)$, so we employ fourth-order compact FDMs in 5 subcases for the case 1 to verify the accuracy and the convergence rate.
		As the remaining 4 subcases only appear when  $q^{(1)}= \pm p^{(1)}$ and $q^{(1)}= \pm \sqrt{3}p^{(1)}$, for the sake of brevity, we only provide  \cref{fig:subcase:2:case:1}--\cref{fig:subcase:4:case:1} which are sufficient for readers to obtain corresponding fourth-order compact FDMs by similarly applying \eqref{f00}--\eqref{Fij:explicit} in the subcase 1.		
	\end{remark}

		\begin{remark}
	Note that transformations in \cref{fig:case:1}--\cref{fig:subcase:4:case:1} are not unique. The first transformation in  \cref{fig:case:1}--\cref{fig:subcase:4:case:1} is the vertical  transformation to avoid the subcase $	q^{(1)}= 0$ to make implementation easy. While second transformations are different to avoid the zero denominator in  the stencil of the proposed FDM.
	In future we plan to establish the optimal transformation to reduce the pollution effect to achieve the smallest error by adopting the truncation error minimization strategy in \cite{FengMichelle2023,FengTrenchea2026}.
	\end{remark}

\begin{figure}[h]
	\begin{subfigure}[b]{0.3\textwidth}
		\hspace{2cm}
		\begin{tikzpicture}[scale = 0.45]
			\node (A) at (0,7) {$u^{(0,0)}$};
			\node (A) at (0,6) {$u^{(0,1)}$};
			\node (A) at (0,5) {$u^{(0,2)}$};
			\node (A) at (0,4) {$u^{(0,3)}$};
			\node (A) at (2,7) {$u^{(1,0)}$};
			\node (A) at (2,6) {$u^{(1,1)}$};
			\node (A) at (2,5) {$u^{(1,2)}$};
			\node (A) at (4,7) {$u^{(2,0)}$};
			\node (A) at (4,6) {$u^{(2,1)}$};
			\node (A) at (6,7) {$u^{(3,0)}$};
			\draw[->, >=stealth, line width=1pt, scale=1.5] (3.2,3.5) -- (7,3.5);
		\end{tikzpicture}
	\end{subfigure}
	\begin{subfigure}[b]{0.3\textwidth}
		\hspace{2cm}
		\begin{tikzpicture}[scale = 0.45]
			\node (A) at (0,7) {$u^{(0,0)}$};
			\node (A) at (0,6) {$u^{(0,1)}$};
			\node (A) at (0,5) {$u^{(0,2)}$};
			\node (A) at (0,4) {$u^{(0,3)}$};
			\node (A) at (2,7) {$u^{(1,0)}$};
			\node (A) at (2,6) {$u^{(1,1)}$};
			\node (A) at (2,5) {$u^{(1,2)}$};
			\draw[->, >=stealth, line width=1pt, scale=1.5] (2.2,3.5) -- (6,3.5);	  	
		\end{tikzpicture}
	\end{subfigure}
	\begin{subfigure}[b]{0.3\textwidth}
		\hspace{1.2cm}
		\vspace{0.4cm}
		\begin{tikzpicture}[scale = 0.45]
			\node (A) at (0,7) {$u^{(1,0)}$};
			\node (A) at (0,6) {$u^{(0,2)}$};
			\node (A) at (0,5) {$u^{(1,2)}$};
		\end{tikzpicture}
	\end{subfigure}	
	\caption
	{The illustration of vertical transformations for the subcase 2 with $q^{(1)}= p^{(1)}$ of the case 1  with $|q^{(1)}/p^{(1)}| \ge 1$. }
	\label{fig:subcase:2:case:1}
\end{figure}		
	\begin{figure}[h]
		 \begin{subfigure}[b]{0.3\textwidth}
			\hspace{2cm}
			\begin{tikzpicture}[scale = 0.45]
				\node (A) at (0,7) {$u^{(0,0)}$};
				\node (A) at (0,6) {$u^{(0,1)}$};
				\node (A) at (0,5) {$u^{(0,2)}$};
				\node (A) at (0,4) {$u^{(0,3)}$};
				\node (A) at (2,7) {$u^{(1,0)}$};
				\node (A) at (2,6) {$u^{(1,1)}$};
				\node (A) at (2,5) {$u^{(1,2)}$};
				\node (A) at (4,7) {$u^{(2,0)}$};
				\node (A) at (4,6) {$u^{(2,1)}$};
				\node (A) at (6,7) {$u^{(3,0)}$};
				\draw[->, >=stealth, line width=1pt, scale=1.5] (3.2,3.5) -- (7,3.5);
			\end{tikzpicture}
		\end{subfigure}
		 \begin{subfigure}[b]{0.3\textwidth}
			\hspace{2cm}
			\begin{tikzpicture}[scale = 0.45]
				\node (A) at (0,7) {$u^{(0,0)}$};
				\node (A) at (0,6) {$u^{(0,1)}$};
				\node (A) at (0,5) {$u^{(0,2)}$};
				\node (A) at (0,4) {$u^{(0,3)}$};
				\node (A) at (2,7) {$u^{(1,0)}$};
				\node (A) at (2,6) {$u^{(1,1)}$};
				\node (A) at (2,5) {$u^{(1,2)}$};
				\draw[->, >=stealth, line width=1pt, scale=1.5] (2.2,3.5) -- (6,3.5);	  	
			\end{tikzpicture}
		\end{subfigure}
		 \begin{subfigure}[b]{0.3\textwidth}
			\hspace{1.2cm}
			\vspace{0.4cm}
			\begin{tikzpicture}[scale = 0.45]
				\node (A) at (0,7) {$u^{(0,1)}$};
				\node (A) at (0,6) {$u^{(0,2)}$};
				\node (A) at (0,5) {$u^{(1,2)}$};
			\end{tikzpicture}
		\end{subfigure}	
		\caption
		{The illustration of vertical transformations for the subcase 3 with $q^{(1)}= -p^{(1)}$ of the case 1 with  $|q^{(1)}/p^{(1)}| \ge 1$. }
		\label{fig:subcase:3:case:1}
	\end{figure}	
\begin{figure}[h]
	\begin{subfigure}[b]{0.3\textwidth}
		\hspace{2cm}
		\begin{tikzpicture}[scale = 0.45]
			\node (A) at (0,7) {$u^{(0,0)}$};
			\node (A) at (0,6) {$u^{(0,1)}$};
			\node (A) at (0,5) {$u^{(0,2)}$};
			\node (A) at (0,4) {$u^{(0,3)}$};
			\node (A) at (2,7) {$u^{(1,0)}$};
			\node (A) at (2,6) {$u^{(1,1)}$};
			\node (A) at (2,5) {$u^{(1,2)}$};
			\node (A) at (4,7) {$u^{(2,0)}$};
			\node (A) at (4,6) {$u^{(2,1)}$};
			\node (A) at (6,7) {$u^{(3,0)}$};
			\draw[->, >=stealth, line width=1pt, scale=1.5] (3.2,3.5) -- (7,3.5);
		\end{tikzpicture}
	\end{subfigure}
	\begin{subfigure}[b]{0.3\textwidth}
		\hspace{2cm}
		\begin{tikzpicture}[scale = 0.45]
			\node (A) at (0,7) {$u^{(0,0)}$};
			\node (A) at (0,6) {$u^{(0,1)}$};
			\node (A) at (0,5) {$u^{(0,2)}$};
			\node (A) at (0,4) {$u^{(0,3)}$};
			\node (A) at (2,7) {$u^{(1,0)}$};
			\node (A) at (2,6) {$u^{(1,1)}$};
			\node (A) at (2,5) {$u^{(1,2)}$};
			\draw[->, >=stealth, line width=1pt, scale=1.5] (2.2,3.5) -- (6,3.5);	  	
		\end{tikzpicture}
	\end{subfigure}
	\begin{subfigure}[b]{0.3\textwidth}
		\hspace{1.2cm}
		\vspace{0.4cm}
		\begin{tikzpicture}[scale = 0.45]
			\node (A) at (0,7) {$u^{(0,1)}$};
			\node (A) at (0,6) {$u^{(1,1)}$};
			\node (A) at (0,5) {$u^{(0,3)}$};
		\end{tikzpicture}
	\end{subfigure}	
	\caption
	{The illustration of vertical transformations for subcases 4 and 5 with $q^{(1)}= \pm \sqrt{3} p^{(1)}$ of the case 1  with $|q^{(1)}/p^{(1)}| \ge 1$. }
	\label{fig:subcase:4:case:1}
\end{figure}

Next, we construct the fourth-order compact FDM at the irregular  stencil center point for the case 2 with $|q^{(1)}/p^{(1)}| \le 1$.

\textbf{Case 2: $|q^{(1)}/p^{(1)}| \le 1$.}
Similar to the case 1, we also need to discuss 5 subcases.
	
\textbf{Subcase 1 of the case 2: $q^{(1)}\neq \pm p^{(1)}$ and $q^{(1)}\neq \pm \frac{\sqrt{3}}{3}p^{(1)}$.} We present the fourth-order compact FDM in the following \cref{thm:FDM:irregular:case:2} directly, while \cref{thm:FDM:irregular:case:2} can be obtained by repeating \eqref{f00}--\eqref{Fij:explicit} with using transformations in \cref{fig:subcase:1:case:2}.
\begin{figure}[h]
	\begin{subfigure}[b]{0.3\textwidth}
		\hspace{2cm}
		\begin{tikzpicture}[scale = 0.45]
			\node (A) at (0,7) {$u^{(0,0)}$};
			\node (A) at (0,6) {$u^{(0,1)}$};
			\node (A) at (0,5) {$u^{(0,2)}$};
			\node (A) at (0,4) {$u^{(0,3)}$};
			\node (A) at (2,7) {$u^{(1,0)}$};
			\node (A) at (2,6) {$u^{(1,1)}$};
			\node (A) at (2,5) {$u^{(1,2)}$};
			\node (A) at (4,7) {$u^{(2,0)}$};
			\node (A) at (4,6) {$u^{(2,1)}$};
			\node (A) at (6,7) {$u^{(3,0)}$};
			\draw[->, >=stealth, line width=1pt, scale=1.5] (3.2,3.5) -- (7,3.5);
		\end{tikzpicture}
	\end{subfigure}
	\begin{subfigure}[b]{0.3\textwidth}
		\hspace{2cm}
		\begin{tikzpicture}[scale = 0.45]
			\node (A) at (0,7) {$u^{(0,0)}$};
			\node (A) at (0,6) {$u^{(1,0)}$};
			\node (A) at (0,5) {$u^{(2,0)}$};
			\node (A) at (0,4) {$u^{(3,0)}$};
			\node (A) at (2,7) {$u^{(0,1)}$};
			\node (A) at (2,6) {$u^{(1,1)}$};
			\node (A) at (2,5) {$u^{(2,1)}$};
			\draw[->, >=stealth, line width=1pt, scale=1.5] (2.2,3.5) -- (6,3.5);	  	
		\end{tikzpicture}
	\end{subfigure}
	\begin{subfigure}[b]{0.3\textwidth}
		\hspace{1.2cm}
		\vspace{0.4cm}
		\begin{tikzpicture}[scale = 0.45]
			\node (A) at (0,7) {$u^{(0,1)}$};
			\node (A) at (0,6) {$u^{(1,1)}$};
			\node (A) at (0,5) {$u^{(2,1)}$};
		\end{tikzpicture}
	\end{subfigure}	
	\caption
	{The illustration of horizontal transformations for the subcase 1 with $q^{(1)}\neq \pm p^{(1)}$ and $q^{(1)}\neq \pm \frac{\sqrt{3}}{3}p^{(1)}$ of the case 2 with $|q^{(1)}/p^{(1)}| \le 1$. }
	\label{fig:subcase:1:case:2}
\end{figure}	
	 \begin{theorem}\label{thm:FDM:irregular:case:2}
	Let $\alpha>0,\vec{\beta},\kappa, \phi,u$ be smooth in $\overline{\Omega}$ in \eqref{model:problem:2D}, $\partial \Omega$ be a smooth boundary curve,  functions $a,b,d,f$ be defined in \eqref{notation:abdf}, and $({x}_i^o,{y}_j^o)\in \Omega$ be the orthogonal projection of the irregular  stencil center point $(x_i,y_j)$. If
	\[
	 ({x}_i^o,{y}_j^o)=(p(t_k^o),q(t_k^o)), \quad |q'(t_k^o)/p'(t_k^o)| \le 1, \quad q'(t_k^o)\ne \pm p'(t_k^o), \quad q'(t_k^o)\ne \pm \frac{\sqrt{3}}{3} p'(t_k^o).
	\]
	Then the following compact FDM  (see  \cref{irregular:FDM:fig} for illustrations)
	\[
	\mathcal{L}_h (u_h)_{i,j} :=\sum_{\substack{r,\ell=-1 \\ (r,\ell)\in S}}^1 C_{r,\ell}\Big|_{(x,y)=(x_i^o,y_j^o)} (u_h)_{i+r,j+\ell}=F_{i,j}\Big|_{(x,y)=(x^o_i,y^o_j)},
	\]
	achieves the fourth-order consistency for the Dirichlet boundary condition	$u =g$ at the point $(x_i^o,y_j^o)\in \partial\Omega$, where the left-hand side  $C_{r,\ell}=\sum_{\sigma=0}^3 c_{r,\ell,\sigma} h^\sigma$ of the stencil  is obtained by  solving 	
	\[
	\sum_{\substack{r,\ell=-1 \\ (r,\ell)\in S}}^1 	A_{r,\ell} c_{r,\ell}=0,\quad
	c_{r,\ell}:=\begin{bmatrix}
		c_{r,\ell,0}\\
		c_{r,\ell,1}\\
		c_{r,\ell,2}\\
		c_{r,\ell,3}
	\end{bmatrix},
	\]
	with the $6 \times 4$ lower triangular matrix $A_{r,\ell}$
	\[
	A_{r,\ell}:=\begin{bmatrix}
\mu_{r}\xi_{3}+\tau_{\ell} & 0 & 0 &  0\\
\mu_{r}^2\xi_{10}+\tau_{\ell}^2\xi_{7}+\mu_{r}\tau_{\ell} & 0 & 0 & 0\\
\mu_{r}\xi_{12}\tau_{\ell}^2-\tau_{\ell}^3/6+\mu_{r}^2\tau_{\ell}/2+\mu_{r}^3\xi_{13} & 0 & 0 & 0\\
\mu_{r}^2\xi_{5}+\tau_{\ell}^2\xi_{1} & \mu_{r}\xi_{3}+\tau_{\ell} & 0 &  0\\
\mu_{r}^3\xi_{11}+\mu_{r}\tau_{\ell}^2\xi_{9}+\tau_{\ell}^3\xi_{8} & \mu_{r}^2\xi_{10}+\tau_{\ell}^2\xi_{7}+\mu_{r}\tau_{\ell} & 0 & 0\\
\mu_{r}^3\xi_{6}+\mu_{r}\tau_{\ell}^2\xi_{4}+\tau_{\ell}^3\xi_{2} & \mu_{r}^2\xi_{5}+\tau_{\ell}^2\xi_{1} & \mu_{r}\xi_{3}+\tau_{\ell} & 0
	\end{bmatrix},
	\]
	\[
	c_{0,0,0}\ne 0,\qquad \sum_{\substack{r,\ell=-1 \\ (r,\ell)\in S}}^1 c_{r, \ell, 0}=1,
	\]
the right-hand side  $F_{i,j}$ of the stencil is defined as 	
\begin{align*}
	F_{i,j}=& \sum_{\substack{r,\ell=-1 \\ (r,\ell)\in S}}^1 \bigg\{ gc_{r, \ell, 0}+\Big\{ c_{r, \ell, 0}g^{(1)}\mu_{r}\omega_{8}+c_{r, \ell, 1}g\Big\}h+\Big\{gc_{r, \ell, 2}+g^{(1)}\omega_{8}\mu_{r}c_{r, \ell, 1}+([\eta_{4}f+g\omega_{4}\\
	 &+g^{(1)}\omega_{10}+g^{(2)}\omega_{15}]\mu_{r}^2+\tau_{\ell}^2[\eta_{1}f+g\omega_{1}+g^{(1)}\omega_{6}+g^{(2)}\omega_{12}])c_{r, \ell, 0}\Big\}h^2+\Big\{([\eta_{5}f+\eta_{7}f^{(0, 1)}\\
	&+\eta_{9}f^{(1, 0)}+g\omega_{5}+g^{(1)}\omega_{11}+g^{(2)}\omega_{16}+g^{(3)}\omega_{18}]\mu_{r}^3+\tau_{\ell}^2(\eta_{3}f+\eta_{6}f^{(0, 1)}+\eta_{8}f^{(1, 0)}+g\omega_{3}\\
	 &+g^{(1)}\omega_{9}+g^{(2)}\omega_{14}+g^{(3)}\omega_{17})\mu_{r}+\tau_{\ell}^3[g\omega_{2}+g^{(1)}\omega_{7}+g^{(2)}\omega_{13}+f\eta_{2}+f^{(0, 1)}/6])c_{r, \ell, 0}\\
	&+c_{r, \ell, 1}(\eta_{4}f+g\omega_{4}+g^{(1)}\omega_{10}+g^{(2)}\omega_{15})\mu_{r}^2+g^{(1)}\omega_{8}\mu_{r}c_{r, \ell, 2}+c_{r, \ell, 1}(\eta_{1}f+g\omega_{1}+g^{(1)}\omega_{6}\\
	 &+g^{(2)}\omega_{12})\tau_{\ell}^2+gc_{r, \ell, 3}\Big\}h^3 \bigg\},
\end{align*}
where 	
	\[
\begin{split}
	& S:=\{ (r,\ell) :   ({x}_i^o+rh,{y}_j^o+\ell h)\in \overline{\Omega} \text{ with } r,\ell=-1,0,1  \},\\
	& w=(x_i-{x}_i^o)/h,\qquad  v=(y_j-{y}_j^o)/h,\qquad \mu_{r}=r+w, \qquad \tau_{\ell}=\ell+v,
\end{split}
\]
\begin{align*}
& \xi_{1} = -(as_{1}+b+s_{2})/2, \\
& \xi_{2} = (b^2+(as_{1}+s_{2})b-s_{1}a^{(0, 1)}-b^{(0, 1)}-d)/6, \\
& \xi_{3} = s_{1},\qquad \xi_{4} = (s_{1}(-d-a^{(1, 0)})-s_{2}a-s_{4}-b^{(1, 0)})/2, \\
& \xi_{5} = s_{2}/2, \qquad \xi_{6} = s_{4}/6, \qquad \xi_{7} = -s_{3}/2, \qquad \xi_{8} = (s_{3}b-a)/6, \\
&  \xi_{9} = -(s_{3}a+b+s_{5})/2, \qquad \xi_{10} = s_{3}/2, \qquad \xi_{11} = s_{5}/6, \\
& \xi_{12} = -s_{6}/2,  \qquad \xi_{13} = s_{6}/6,\\
& \eta_{1} = (z_{3}+1)/2,\qquad \eta_{2} = -(z_{3}+1)b/6, \qquad \eta_{3} = (z_{3}a+z_{7})/2, \\
&  \eta_{4} = -z_{3}/2,   \qquad \eta_{5} = -z_{7}/6, \qquad \eta_{6} = z_{8}/2, \qquad \eta_{7} = -z_{8}/6, \\
&  \eta_{8} = (z_{10}+1)/2, \qquad \eta_{9} = -z_{10}/6,\\
& \omega_{1} = (z_{2}-d)/2,\qquad  \omega_{2} = ([-z_{2}+d]b-d^{(0, 1)})/6,    \\
& \omega_{3} = (z_{2}a-d^{(1, 0)}+z_{6})/2, \qquad \omega_{4} = -z_{2}/2,\qquad \omega_{5} = -z_{6}/6, \\
& \omega_{6} = (az_{1}+z_{4})/2, \qquad \omega_{7} = (a^{(0, 1)}z_{1}-[az_{1}+z_{4}]b)/6,  \\
&  \omega_{8} = -z_{1}, \qquad \omega_{9} = ([a^{(1, 0)}+d]z_{1}+z_{4}a+z_{9})/2, \\
&  \omega_{10} = -z_{4}/2,  \qquad  \omega_{11} = -z_{9}/6,  \qquad \omega_{12} = z_{5}/2,\\
&   \omega_{13} = -z_{5}b/6,  \qquad  \omega_{14} = (z_{5}a+z_{11})/2,  \qquad  \omega_{15} = -z_{5}/2,   \\
&   \omega_{16} = -z_{11}/6,  \qquad \omega_{17} = z_{12}/2, \qquad \omega_{18} = -z_{12}/6,
\end{align*}
with
\begin{align*}
	s_1 =& -q^{(1)}/p^{(1)}, \\
	s_2 =& \big\{(b[q^{(1)}]^2-q^{(2)})p^{(1)}-a[q^{(1)}]^3\\
	 &+p^{(2)}q^{(1)}\big\}/\big\{p^{(1)}([p^{(1)}]^2-[q^{(1)}]^2)\big\}, \\
	\quad s_3 =& -2p^{(1)}q^{(1)}/\big\{[p^{(1)}]^2-[q^{(1)}]^2\big\}, \\
	\quad s_4 = &\big\{a^{(0, 1)}[q^{(1)}]^6-3(a^2-a^{(1, 0)}+b^{(0, 1)}/3-2d/3)p^{(1)}[q^{(1)}]^5\\
	&+([p^{(1)}]^2(4ab-a^{(0, 1)}-3b^{(1, 0)})-bp^{(2)})[q^{(1)}]^4+([p^{(1)}]^3(b^{(0, 1)}\\
	&-b^2-3a^{(1, 0)}-2d)+(6ap^{(2)}+bq^{(2)})p^{(1)}-p^{(3)})[q^{(1)}]^3\\
	&+(3b^{(1, 0)}[p^{(1)}]^4-3(2aq^{(2)}+bp^{(2)})[p^{(1)}]^2+p^{(1)}q^{(3)}\\
	 &+3p^{(2)}q^{(2)})[q^{(1)}]^2+3p^{(1)}(b[p^{(1)}]^2q^{(2)}+p^{(1)}p^{(3)}/3\\
	 &-[p^{(2)}]^2-[q^{(2)}]^2)q^{(1)}-[p^{(1)}]^3q^{(3)}+3[p^{(1)}]^2p^{(2)}q^{(2)}\big\}\\
	 &/\big\{[p^{(1)}]^2([p^{(1)}]^4-4[p^{(1)}]^2[q^{(1)}]^2+3[q^{(1)}]^4)\big\}, \\
	s_5 =&  \big\{3b[p^{(1)}]^3[q^{(1)}]^2-5a[p^{(1)}]^2[q^{(1)}]^3-a[q^{(1)}]^5-bp^{(1)}[q^{(1)}]^4\\
	 &-3[p^{(1)}]^3q^{(2)}+3[p^{(1)}]^2p^{(2)}q^{(1)}-3p^{(1)}[q^{(1)}]^2q^{(2)}+3p^{(2)}[q^{(1)}]^3\big\}\\
	 &/\big\{p^{(1)}([p^{(1)}]^4-4[p^{(1)}]^2[q^{(1)}]^2+3[q^{(1)}]^4)\big\},\\
	s_6 = & \big\{[q^{(1)}]^3-3[p^{(1)}]^2q^{(1)}\big\}/\big\{p^{(1)}([p^{(1)}]^2-3[q^{(1)}]^2)\big\},
\end{align*}
and
\begin{align*}
	z_{1} =& -1/p^{(1)}, \\
	z_2 =& -d[q^{(1)}]^2/\big\{[p^{(1)}]^2-[q^{(1)}]^2\big\}, \\
	z_3 =& [q^{(1)}]^2/\big\{[p^{(1)}]^2-[q^{(1)}]^2\big\},\\
	z_4 =& \big\{-a[q^{(1)}]^2+p^{(2)}\big\}/\big\{p^{(1)}([p^{(1)}]^2-[q^{(1)}]^2)\big\},\\
	z_5 =& -1/\big\{[p^{(1)}]^2-[q^{(1)}]^2\big\}, \\
	z_6 =& -q^{(1)}\big\{ 3adp^{(1)}[q^{(1)}]^3-bd[p^{(1)}]^2[q^{(1)}]^2\\
	&+d^{(0, 1)}[p^{(1)}]^2[q^{(1)}]^2-d^{(0, 1)}[q^{(1)}]^4+3d^{(1, 0)}[p^{(1)}]^3q^{(1)}\\
	&-3d^{(1, 0)}p^{(1)}[q^{(1)}]^3+3d[p^{(1)}]^2q^{(2)}-3dp^{(1)}p^{(2)}q^{(1)}\big\}\\
	 &/\big\{p^{(1)}([p^{(1)}]^4-4[p^{(1)}]^2[q^{(1)}]^2+3[q^{(1)}]^4)\big\}, \\
	z_7=& 3q^{(1)}\big\{a[q^{(1)}]^3-bp^{(1)}[q^{(1)}]^2/3+p^{(1)}q^{(2)}-p^{(2)}q^{(1)}\big\}\\
	 &/\big\{[p^{(1)}]^4-4[p^{(1)}]^2[q^{(1)}]^2+3[q^{(1)}]^4\big\},\\
	z_8 = & [q^{(1)}]^3/\big\{p^{(1)}([p^{(1)}]^2-3[q^{(1)}]^2)\big\},\\
	z_9= & \big\{a^{(0, 1)}[q^{(1)}]^5-3p^{(1)}(a^2-a^{(1, 0)}-d)[q^{(1)}]^4\\
	&+((ab-a^{(0, 1)})[p^{(1)}]^2-bp^{(2)})[q^{(1)}]^3+(-3(a^{(1, 0)}\\
	 &+d)[p^{(1)}]^3+6ap^{(1)}p^{(2)}-p^{(3)})[q^{(1)}]^2-3q^{(2)}(a[p^{(1)}]^2\\
	 &-p^{(2)})q^{(1)}+[p^{(1)}]^2p^{(3)}-3p^{(1)}[p^{(2)}]^2\big\}/\big\{[p^{(1)}]^2([p^{(1)}]^4\\
	 &-4[p^{(1)}]^2[q^{(1)}]^2+3[q^{(1)}]^4)\big\},  \\
	z_{10} = & 3[q^{(1)}]^2/\big\{[p^{(1)}]^2-3[q^{(1)}]^2\big\}, \\
	z_{11} =& \big\{b[q^{(1)}]^3-3ap^{(1)}[q^{(1)}]^2+3p^{(1)}p^{(2)}-3q^{(1)}q^{(2)}\big\}\\
	 &/\big\{p^{(1)}([p^{(1)}]^4-4[p^{(1)}]^2[q^{(1)}]^2+3[q^{(1)}]^4)\big\}, \\
	z_{12} =& -1/\big\{p^{(1)}([p^{(1)}]^2-3[q^{(1)}]^2)\big\}.
\end{align*}
\end{theorem}
Similar to subcases 2 to 5 of the case 1, we also provide \cref{fig:subcase:2:case:2}--\cref{fig:subcase:4:case:2} for subcases 2 to 5 of the case 2 which are sufficient to derive corresponding fourth-order compact FDMs at irregular  stencil center points in the remaining 4 subcases by adopting \eqref{f00}--\eqref{Fij:explicit} similarly.
\begin{figure}[h]
	\begin{subfigure}[b]{0.3\textwidth}
		\hspace{2cm}
		\begin{tikzpicture}[scale = 0.45]
			\node (A) at (0,7) {$u^{(0,0)}$};
			\node (A) at (0,6) {$u^{(0,1)}$};
			\node (A) at (0,5) {$u^{(0,2)}$};
			\node (A) at (0,4) {$u^{(0,3)}$};
			\node (A) at (2,7) {$u^{(1,0)}$};
			\node (A) at (2,6) {$u^{(1,1)}$};
			\node (A) at (2,5) {$u^{(1,2)}$};
			\node (A) at (4,7) {$u^{(2,0)}$};
			\node (A) at (4,6) {$u^{(2,1)}$};
			\node (A) at (6,7) {$u^{(3,0)}$};
			\draw[->, >=stealth, line width=1pt, scale=1.5] (3.2,3.5) -- (7,3.5);
		\end{tikzpicture}
	\end{subfigure}
	\begin{subfigure}[b]{0.3\textwidth}
		\hspace{2cm}
		\begin{tikzpicture}[scale = 0.45]
			\node (A) at (0,7) {$u^{(0,0)}$};
			\node (A) at (0,6) {$u^{(1,0)}$};
			\node (A) at (0,5) {$u^{(2,0)}$};
			\node (A) at (0,4) {$u^{(3,0)}$};
			\node (A) at (2,7) {$u^{(0,1)}$};
			\node (A) at (2,6) {$u^{(1,1)}$};
			\node (A) at (2,5) {$u^{(2,1)}$};
			\draw[->, >=stealth, line width=1pt, scale=1.5] (2.2,3.5) -- (6,3.5);	  	
		\end{tikzpicture}
	\end{subfigure}
	\begin{subfigure}[b]{0.3\textwidth}
		\hspace{1.2cm}
		\vspace{0.4cm}
		\begin{tikzpicture}[scale = 0.45]
			\node (A) at (0,7) {$u^{(1,0)}$};
			\node (A) at (0,6) {$u^{(2,0)}$};
			\node (A) at (0,5) {$u^{(3,0)}$};
		\end{tikzpicture}
	\end{subfigure}	
	\caption
	{The illustration of horizontal transformations for the subcase 2 with $q^{(1)}= p^{(1)}$ of the case 2 with $|q^{(1)}/p^{(1)}| \le 1$. }
	\label{fig:subcase:2:case:2}
\end{figure}	
\begin{figure}[h]
	\begin{subfigure}[b]{0.3\textwidth}
		\hspace{2cm}
		\begin{tikzpicture}[scale = 0.45]
			\node (A) at (0,7) {$u^{(0,0)}$};
			\node (A) at (0,6) {$u^{(0,1)}$};
			\node (A) at (0,5) {$u^{(0,2)}$};
			\node (A) at (0,4) {$u^{(0,3)}$};
			\node (A) at (2,7) {$u^{(1,0)}$};
			\node (A) at (2,6) {$u^{(1,1)}$};
			\node (A) at (2,5) {$u^{(1,2)}$};
			\node (A) at (4,7) {$u^{(2,0)}$};
			\node (A) at (4,6) {$u^{(2,1)}$};
			\node (A) at (6,7) {$u^{(3,0)}$};
			\draw[->, >=stealth, line width=1pt, scale=1.5] (3.2,3.5) -- (7,3.5);
		\end{tikzpicture}
	\end{subfigure}
	\begin{subfigure}[b]{0.3\textwidth}
		\hspace{2cm}
		\begin{tikzpicture}[scale = 0.45]
			\node (A) at (0,7) {$u^{(0,0)}$};
			\node (A) at (0,6) {$u^{(1,0)}$};
			\node (A) at (0,5) {$u^{(2,0)}$};
			\node (A) at (0,4) {$u^{(3,0)}$};
			\node (A) at (2,7) {$u^{(0,1)}$};
			\node (A) at (2,6) {$u^{(1,1)}$};
			\node (A) at (2,5) {$u^{(2,1)}$};
			\draw[->, >=stealth, line width=1pt, scale=1.5] (2.2,3.5) -- (6,3.5);	  	
		\end{tikzpicture}
	\end{subfigure}
	\begin{subfigure}[b]{0.3\textwidth}
		\hspace{1.2cm}
		\vspace{0.4cm}
		\begin{tikzpicture}[scale = 0.45]
			\node (A) at (0,7) {$u^{(0,1)}$};
			\node (A) at (0,6) {$u^{(2,0)}$};
			\node (A) at (0,5) {$u^{(3,0)}$};
		\end{tikzpicture}
	\end{subfigure}	
	\caption
	{The illustration of horizontal transformations for the subcase 3 with  $q^{(1)}= - p^{(1)}$ of the case 2  with $|q^{(1)}/p^{(1)}| \le 1$. }
	\label{fig:subcase:3:case:2}
\end{figure}	
\begin{figure}[h]
	\begin{subfigure}[b]{0.3\textwidth}
		\hspace{2cm}
		\begin{tikzpicture}[scale = 0.45]
			\node (A) at (0,7) {$u^{(0,0)}$};
			\node (A) at (0,6) {$u^{(0,1)}$};
			\node (A) at (0,5) {$u^{(0,2)}$};
			\node (A) at (0,4) {$u^{(0,3)}$};
			\node (A) at (2,7) {$u^{(1,0)}$};
			\node (A) at (2,6) {$u^{(1,1)}$};
			\node (A) at (2,5) {$u^{(1,2)}$};
			\node (A) at (4,7) {$u^{(2,0)}$};
			\node (A) at (4,6) {$u^{(2,1)}$};
			\node (A) at (6,7) {$u^{(3,0)}$};
			\draw[->, >=stealth, line width=1pt, scale=1.5] (3.2,3.5) -- (7,3.5);
		\end{tikzpicture}
	\end{subfigure}
	\begin{subfigure}[b]{0.3\textwidth}
		\hspace{2cm}
		\begin{tikzpicture}[scale = 0.45]
			\node (A) at (0,7) {$u^{(0,0)}$};
			\node (A) at (0,6) {$u^{(1,0)}$};
			\node (A) at (0,5) {$u^{(2,0)}$};
			\node (A) at (0,4) {$u^{(3,0)}$};
			\node (A) at (2,7) {$u^{(0,1)}$};
			\node (A) at (2,6) {$u^{(1,1)}$};
			\node (A) at (2,5) {$u^{(2,1)}$};
			\draw[->, >=stealth, line width=1pt, scale=1.5] (2.2,3.5) -- (6,3.5);	  	
		\end{tikzpicture}
	\end{subfigure}
	\begin{subfigure}[b]{0.3\textwidth}
		\hspace{1.2cm}
		\vspace{0.4cm}
		\begin{tikzpicture}[scale = 0.45]
			\node (A) at (0,7) {$u^{(0,1)}$};
			\node (A) at (0,6) {$u^{(1,1)}$};
			\node (A) at (0,5) {$u^{(3,0)}$};
		\end{tikzpicture}
	\end{subfigure}	
	\caption
	{The illustration of horizontal transformations for subcases 4 and 5 with $q^{(1)}= \pm \frac{\sqrt{3}}{3}p^{(1)}$  of the case 2  with $|q^{(1)}/p^{(1)}| \le 1$. }
	\label{fig:subcase:4:case:2}
\end{figure}

	\section{Numerical experiments}\label{sec:Numeri}

	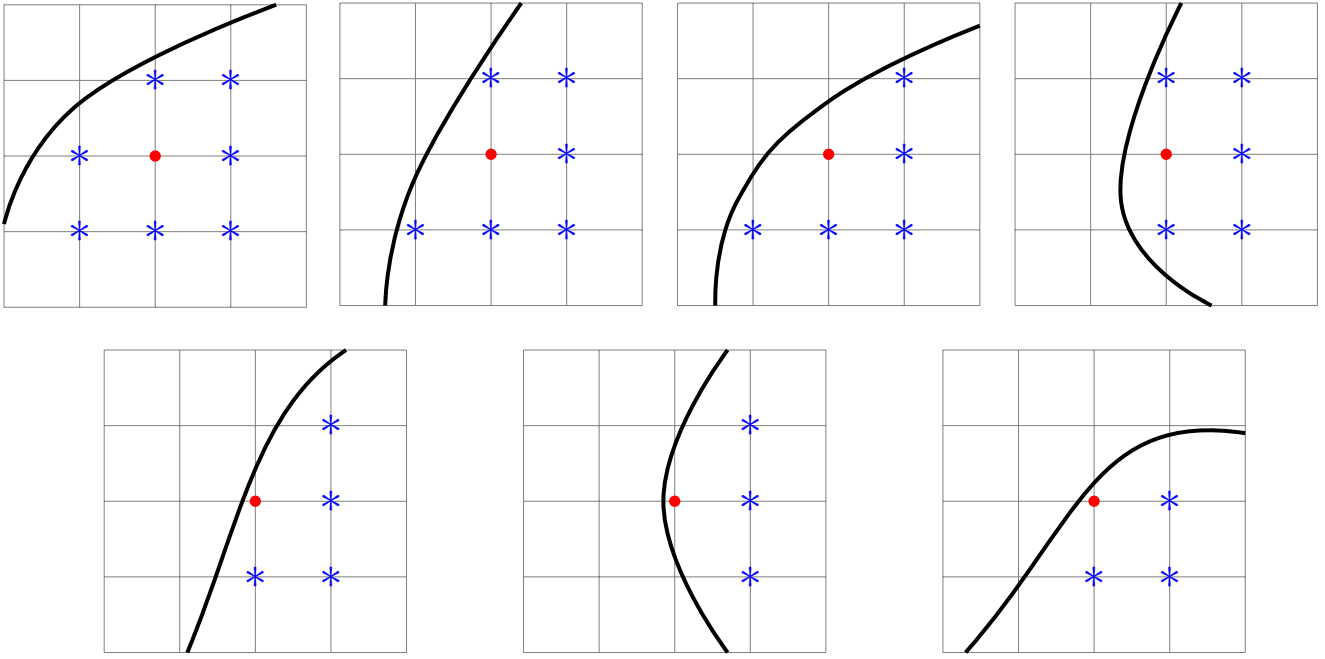
\begin{figure}[h]
	\centering
	\begin{subfigure}[b]{0.24\textwidth}
		\begin{tikzpicture}[scale = 1]
			\draw[help lines,step = 1]
			(0,0) grid (4,4);
			\node at (1,1) { \Large \color{blue}$*$};
			\node at (1,2) { \Large \color{blue}$*$};
			\node at (2,1) { \Large \color{blue}$*$};
			\node at (2,2)[circle,fill,inner sep=1.5pt,color=red]{};
			\node at (2,3){ \Large \color{blue}$*$};
			\node at (3,1) { \Large \color{blue}$*$};
			\node at (3,2) { \Large \color{blue}$*$};
			\node at (3,3) { \Large \color{blue}$*$};
			\draw[line width=1.5pt, black]  plot [smooth,tension=0.8]
			coordinates {(0,1.1) (1,2.7) (3.6,4)};
		\end{tikzpicture}
	\end{subfigure}
	\begin{subfigure}[b]{0.24\textwidth}
		\begin{tikzpicture}[scale = 1]
			\draw[help lines,step = 1]
			(0,0) grid (4,4);
			\node at (1,1) { \Large \color{blue}$*$};
			\node at (2,1) { \Large \color{blue}$*$};
			\node at (2,2)[circle,fill,inner sep=1.5pt,color=red]{};
			\node at (2,3){ \Large \color{blue}$*$};
			\node at (3,1) { \Large \color{blue}$*$};
			\node at (3,2) { \Large \color{blue}$*$};
			\node at (3,3) { \Large \color{blue}$*$};
			\draw[line width=1.5pt, black]  plot [smooth,tension=0.8]
			coordinates {(0.6,0) (1,1.7) (2.4,4)};
		\end{tikzpicture}
	\end{subfigure}
	\begin{subfigure}[b]{0.24\textwidth}
		\begin{tikzpicture}[scale = 1]
			\draw[help lines,step = 1]
			(0,0) grid (4,4);
			\node at (1,1) { \Large \color{blue}$*$};
			\node at (2,1) { \Large \color{blue}$*$};
			\node at (2,2)[circle,fill,inner sep=1.5pt,color=red]{};
			\node at (3,1) { \Large \color{blue}$*$};
			\node at (3,2) { \Large \color{blue}$*$};
			\node at (3,3) { \Large \color{blue}$*$};
			\draw[line width=1.5pt, black]  plot [smooth,tension=0.8]
			coordinates {(0.5,0) (0.8,1.4) (2,2.7) (4,3.7)};
		\end{tikzpicture}
	\end{subfigure}
	\begin{subfigure}[b]{0.24\textwidth}
		\begin{tikzpicture}[scale = 1]
			\draw[help lines,step = 1]
			(0,0) grid (4,4);
			\node at (2,1) { \Large \color{blue}$*$};
			\node at (2,2)[circle,fill,inner sep=1.5pt,color=red]{};
			\node at (2,3){ \Large \color{blue}$*$};
			\node at (3,1) { \Large \color{blue}$*$};
			\node at (3,2) { \Large \color{blue}$*$};
			\node at (3,3) { \Large \color{blue}$*$};
			\draw[line width=1.5pt, black]  plot [smooth,tension=0.8]
			coordinates {(2.6,0) (1.4,1.4) (2.2,4)};
		\end{tikzpicture}
	\end{subfigure}	
	\begin{subfigure}[b]{0.3\textwidth}
		\vspace{0.5cm} \hspace{0.6cm}
		\begin{tikzpicture}[scale = 1]
			\draw[help lines,step = 1]
			(0,0) grid (4,4);
			\node at (2,1) { \Large \color{blue}$*$};
			\node at (2,2)[circle,fill,inner sep=1.5pt,color=red]{};
			\node at (3,1) { \Large \color{blue}$*$};
			\node at (3,2) { \Large \color{blue}$*$};
			\node at (3,3) { \Large \color{blue}$*$};
			\draw[line width=1.5pt, black] (1.1,0) .. controls (1.8,1.7) and (2,3.2) .. (3.2,4);
		\end{tikzpicture}
	\end{subfigure}
		\begin{subfigure}[b]{0.3\textwidth}
		\vspace{0.5cm} \hspace{0.6cm}
		\begin{tikzpicture}[scale = 1]
	\draw[help lines,step = 1]
	(0,0) grid (4,4);
	\node at (2,2)[circle,fill,inner sep=1.5pt,color=red]{};
	\node at (3,1) { \Large \color{blue}$*$};
	\node at (3,2) { \Large \color{blue}$*$};
	\node at (3,3) { \Large \color{blue}$*$};
	\draw[line width=1.5pt, black]  plot [smooth,tension=0.8]
coordinates {(2.7,0) (1.85,2) (2.7,4)};
\end{tikzpicture}
	\end{subfigure}
\begin{subfigure}[b]{0.3\textwidth}
	\vspace{0.5cm} \hspace{0.6cm}
	\begin{tikzpicture}[scale = 1]
		\draw[help lines,step = 1]
		(0,0) grid (4,4);
		\node at (2,1) { \Large \color{blue}$*$};
		\node at (2,2)[circle,fill,inner sep=1.5pt,color=red]{};
		\node at (3,1) { \Large \color{blue}$*$};
		\node at (3,2) { \Large \color{blue}$*$};
		\draw[line width=1.5pt, black] (0.3,0) .. controls (1.8,1.7) and (2,3.2) .. (4,2.9);
	\end{tikzpicture}
\end{subfigure}	
	\caption
	{Up to symmetric
		transformations and rotations, compact FDMs at irregular  stencil center points only have the above 7 configurations for any complex boundary curves  when the mesh size $h$ is reasonably small, where red points are irregular  stencil center points,  blue star points are remaining grid points used in compact FDMs, and black curves are boundaries. See \cite[Figs. 3-5]{Feng2022} for more configurations if $h$ is not sufficiently small.
}
	\label{5:configurations:fig}
\end{figure}	
	
	Recall that we define that $u_{i,j}=u(x_i,y_j)$, and $(u_h)_{i,j}=$ the value of  $u_h$ at the grid point $(x_i,y_j)$, where $\Omega$ is discretized  by the following uniform Cartesian grid:
	\[
	\overline{\Omega}_h:=\big\{(x_i,y_j) \cap \overline{\Omega}\quad : \quad (x_i,y_j):=(i h,jh), \quad i,j=0,\pm 1, \pm 2 \ldots,  \quad \text{and} \quad h:=1/N, \quad N\in \N\big\}.
	\]
Then we define the following $l_2$ and $l_{\infty}$ norms of errors to verify the convergence rate of the proposed fourth-order compact FDM	
	\[
	\begin{split}
		& \frac{\|u_h-u\|_2}{\|u\|_2}		 :=\sqrt{\frac{\sum_{(x_i,y_j)\in \overline{\Omega}_h} \left|(u_h)_{i,j}-u(x_i,y_j)\right|^2}{\sum_{(x_i,y_j)\in \overline{\Omega}_h} \left|u(x_i,y_j)\right|^2}},\\
		& \|u_h-u\|_\infty
		:=\max_{(x_i,y_j)\in \overline{\Omega}_h} \left|(u_h)_{i,j}-u(x_i,y_j)\right|.
	\end{split}
	\]

	In the following \cref{Example:1,Example:2,Example:3}, the domain $\Omega$ is very thin and/or the boundary curve $\partial\Omega$ is sharply varying and/or exhibits the high frequency.

	\begin{figure}[htbp]
		\centering
		 \begin{subfigure}[b]{0.32\textwidth}
			 \includegraphics[width=5.7cm,height=5.7cm]{curve1.pdf}
		\end{subfigure}
		 \begin{subfigure}[b]{0.32\textwidth}
			 \includegraphics[width=6cm,height=6cm]{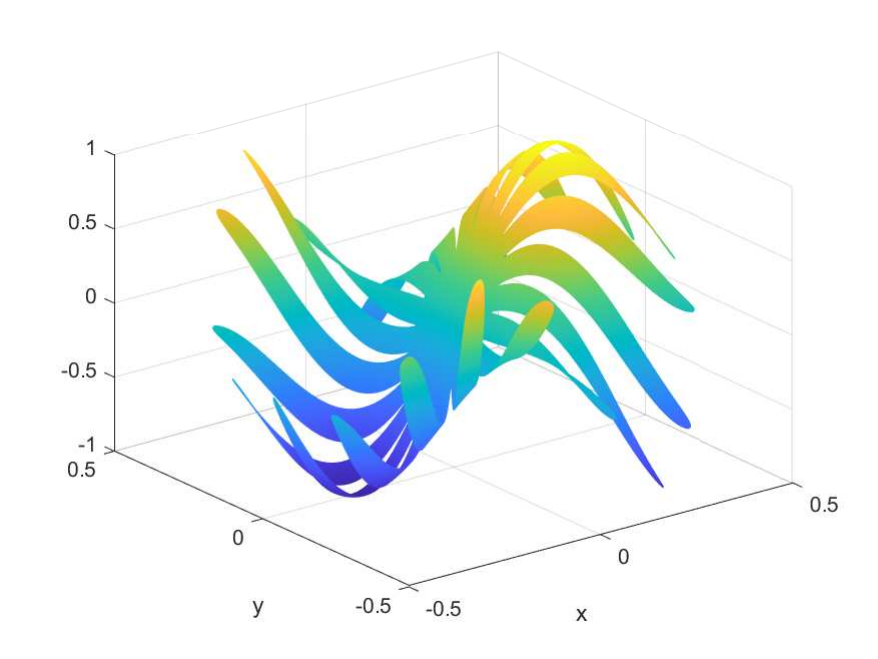}
		\end{subfigure}
		 \begin{subfigure}[b]{0.32\textwidth}
			 \includegraphics[width=6cm,height=6cm]{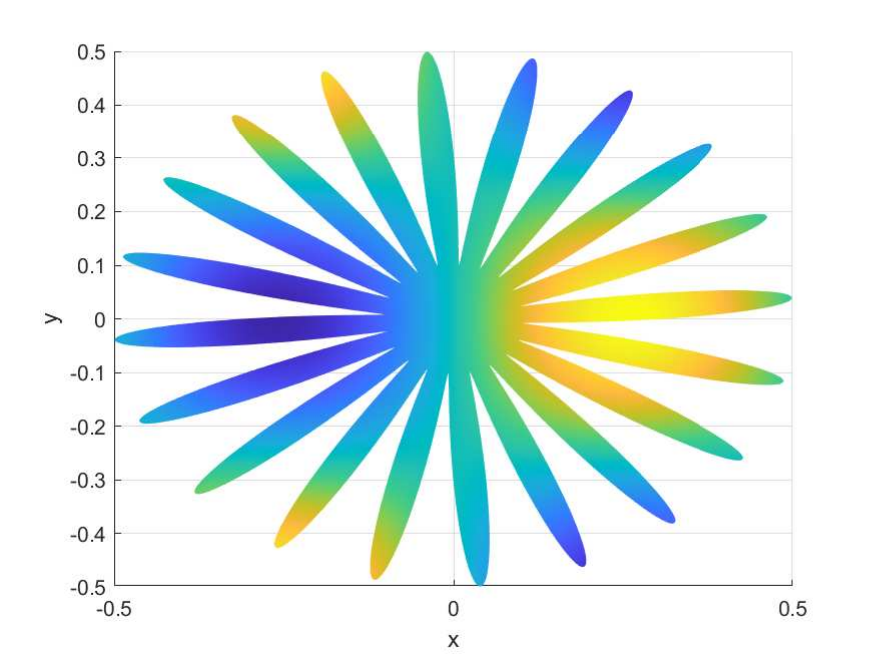}
		\end{subfigure}
		 \begin{subfigure}[b]{0.37\textwidth}
			 \includegraphics[width=6cm,height=6cm]{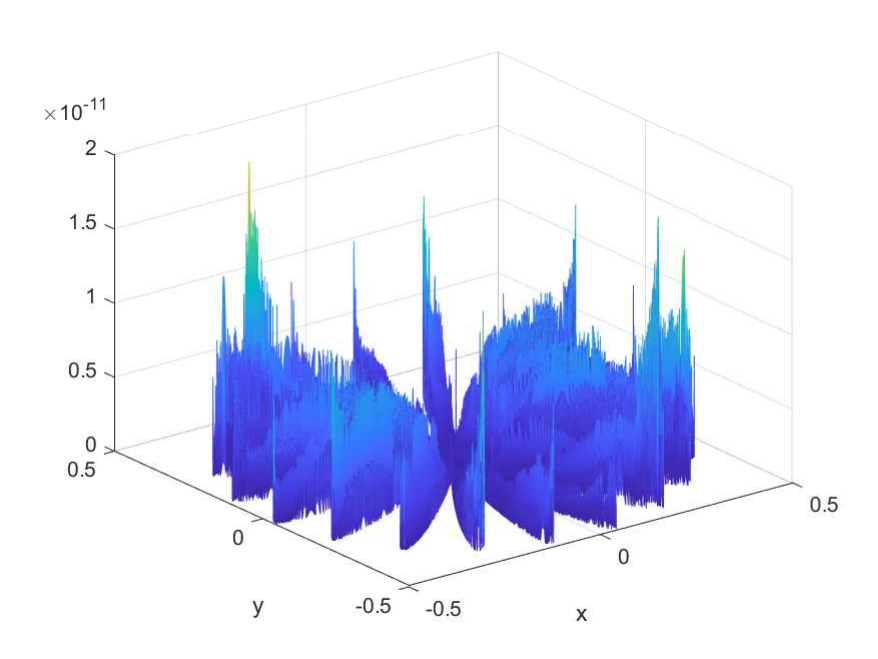}
		\end{subfigure}
		 \begin{subfigure}[b]{0.32\textwidth}
			 \includegraphics[width=6cm,height=6cm]{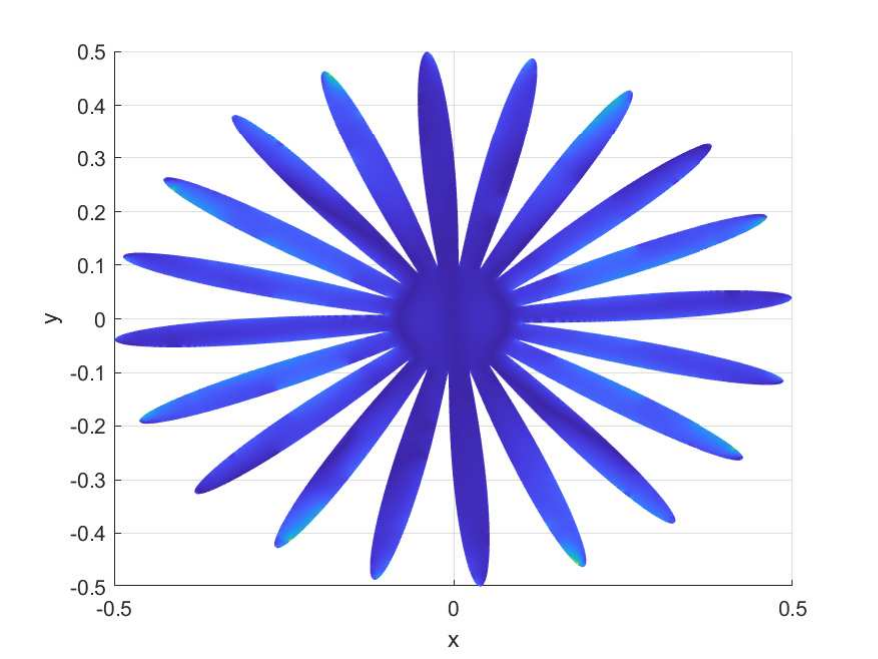}
		\end{subfigure}
		\caption
		{Performance in \cref{Example:1} of the proposed fourth-order compact FDM. The domain $\Omega$ that is enclosed by a 20-leaf boundary curve $\partial \Omega$ (left panel in the first row), $u_h$  on $\overline{\Omega}_h$ with $h=1/2^{12}$ (middle and right panels in the first row), and $|u_h-u|$   on $\overline{\Omega}_h$  with $h=1/2^{12}$ (left and right panels in the second row).}
		\label{Example:1:fig}
	\end{figure}

	\begin{example}\label{Example:1}
		\normalfont
		The functions  of the convection-diffusion-reaction equation \eqref{model:problem:2D} on a curved domain $\Omega$ are given by
		\begin{align*}
			&u=\sin(6 x)\cos(6 y),\qquad \alpha= \exp(x+y) ,\qquad  \\
			&\beta_1= \exp(x-y), \qquad \beta_2= \cos(x)\cos(y), \qquad \kappa=\exp(x+y), \\
			& \Omega=\{(x(t), y(t)) \ : \ x(t)^2+y(t)^2< [0.3+0.2\sin(20t)]^2, \quad t\in [0,2\pi)\},
		\end{align*}
		the source term $\phi$	and the Dirichlet boundary function $g$  are obtained by plugging above functions into \eqref{model:problem:2D}. The numerical results are presented in \cref{Example:1:table} and \cref{Example:1:fig}.	
	\end{example}
	
	\begin{table}[htbp]
		\caption{Performance in \cref{Example:1} of the proposed fourth-order compact FDM.}
		\centering
		 {\renewcommand{\arraystretch}{1.0}
			\scalebox{1}{
				 \setlength{\tabcolsep}{5mm}{
					 \begin{tabular}{c|c|c|c|c}
						\hline
						$h$ &   $\frac{\|u_h-u\|_2}{\|u\|_2}$    &order &   	   $\|u_{h}-u\|_{\infty}$    &order \\
						\hline
$1/2^5$  &  9.0345E-01  &    &  4.4766E+00  &  \\
$1/2^6$  &  5.4344E-02  &  4.06  &  5.9199E-01  &  2.92\\
$1/2^7$  &  6.8879E-04  &  6.30  &  3.9511E-03  &  7.23\\
$1/2^8$  &  3.5861E-07  &  10.91  &  2.2467E-06  &  10.78\\
$1/2^9$  &  2.2640E-08  &  3.99  &  8.0102E-08  &  4.81\\
$1/2^{10}$  &  1.3953E-09  &  4.02  &  5.1896E-09  &  3.95\\
$1/2^{11}$  &  8.7934E-11  &  3.99  &  3.3360E-10  &  3.96\\
$1/2^{12}$  &  4.8243E-12  &  4.19  &  1.9917E-11  &  4.07\\
						\hline
		\end{tabular}}}}
		\label{Example:1:table}
	\end{table}

	\begin{remark}
	Up to symmetric
	transformations and rotations, the stencil of the proposed fourth-order compact FDM at the irregular  stencil center point only has 7 configurations depicted in \cref{5:configurations:fig} for any complicated boundary curves, when the mesh size $h$ is reasonably small. So  we notice that the convergence rate increases to 11 when $h=1/2^{8}$, and then the convergence rate is stabilized at 4 when $h< 1/2^{8}$ (the reasonably fine meshes) in \cref{Example:1:table} for \cref{Example:1}. Similar observations can be found in \cref{Example:2:table:1,Example:2:table:2} for \cref{Example:2} and \cref{Example:3:table:1,Example:3:table:2,Example:3:table:3,Example:3:table:4} for \cref{Example:3}.
	\end{remark}

	\begin{figure}[htbp]
		\centering
		 \begin{subfigure}[b]{0.32\textwidth}
			 \includegraphics[width=5.7cm,height=5.7cm]{curve21.pdf}
		\end{subfigure}
		 \begin{subfigure}[b]{0.32\textwidth}
			 \includegraphics[width=6cm,height=6cm]{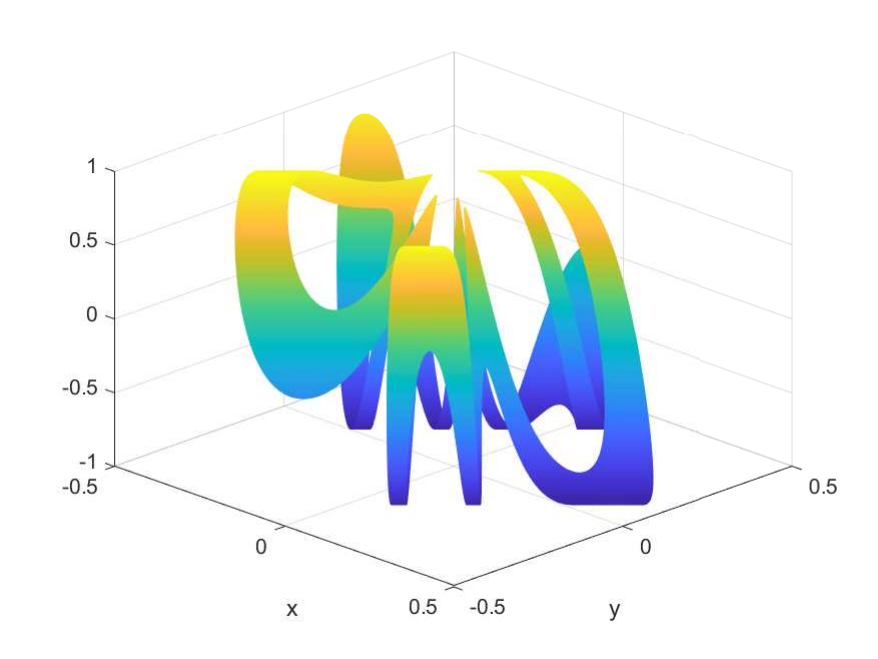}
		\end{subfigure}
		 \begin{subfigure}[b]{0.32\textwidth}
			 \includegraphics[width=6cm,height=6cm]{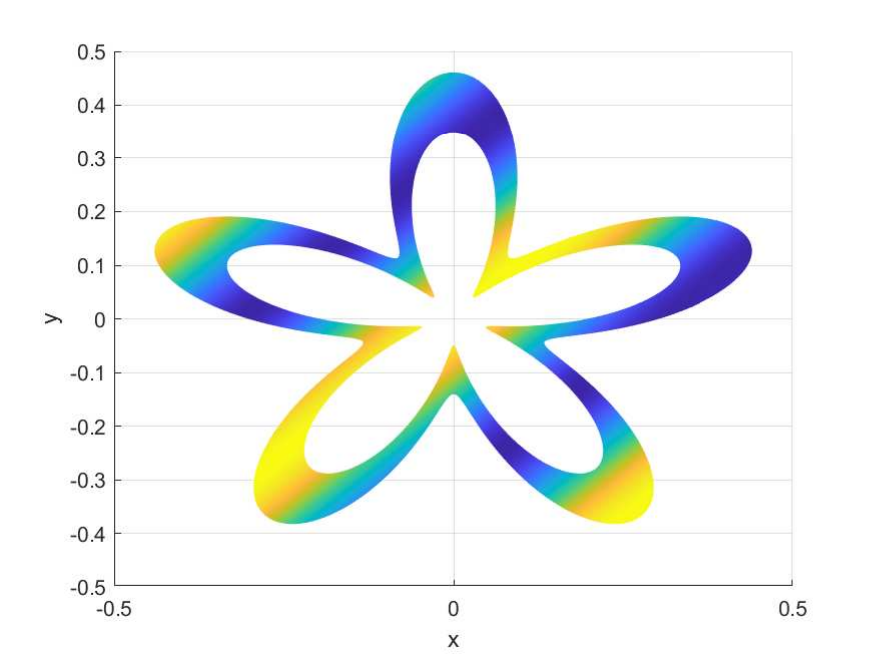}
		\end{subfigure}
		 \begin{subfigure}[b]{0.37\textwidth}
			 \includegraphics[width=6cm,height=6cm]{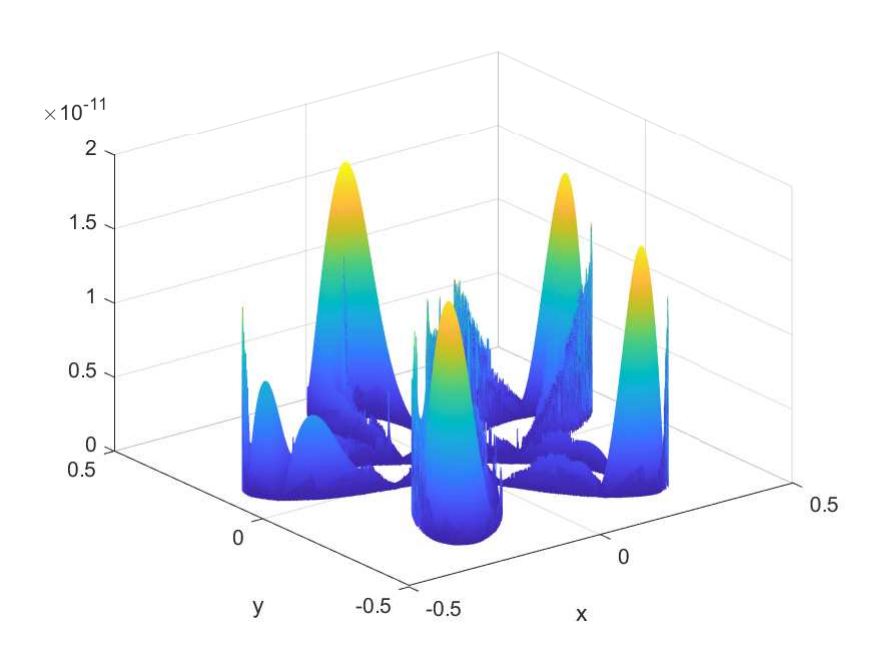}
		\end{subfigure}
		 \begin{subfigure}[b]{0.32\textwidth}
			 \includegraphics[width=6cm,height=6cm]{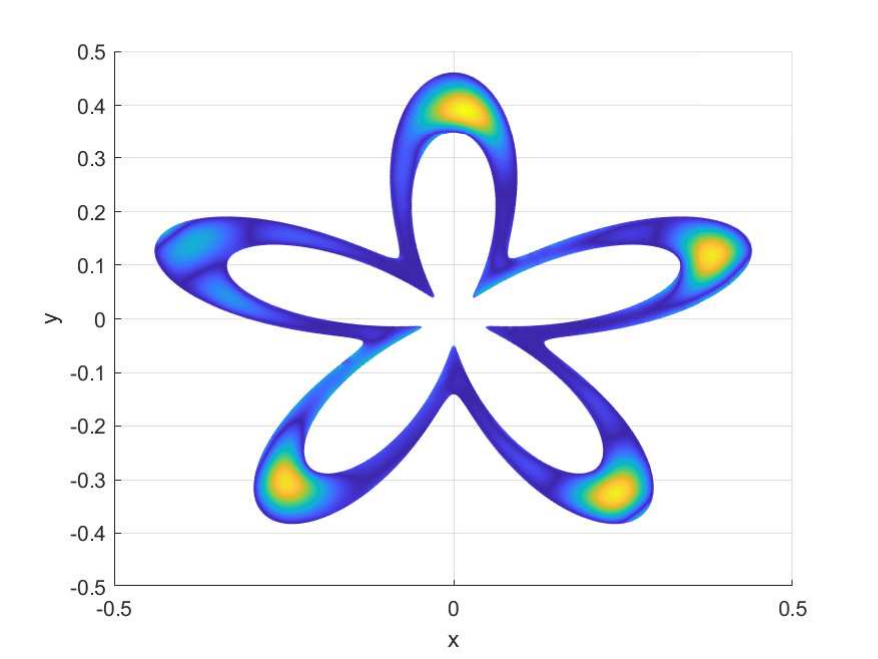}
		\end{subfigure}
		\caption
		{Performance in \cref{Example:2} of the proposed fourth-order compact FDM with $u=\cos(10(x-y))$. The domain $\Omega=\{(x(t), y(t)) : [0.2+0.15\sin(5t)]^2 < x(t)^2+y(t)^2 < [0.3+0.16\sin(5t)]^2, \ t\in [0,2\pi)\}$ that is enclosed by two 5-leaf boundary curves $\partial \Omega$ (left panel in the first row), $u_h$  on $\overline{\Omega}_h$ with $h=1/2^{13}$ (middle and right panels in the first row), and $|u_h-u|$   on $\overline{\Omega}_h$  with $h=1/2^{13}$ (left and right panels in the second row).}
		\label{Example:2:fig:1}
	\end{figure}

	\begin{figure}[htbp]
		\centering
		 \begin{subfigure}[b]{0.32\textwidth}
			 \includegraphics[width=5.7cm,height=5.7cm]{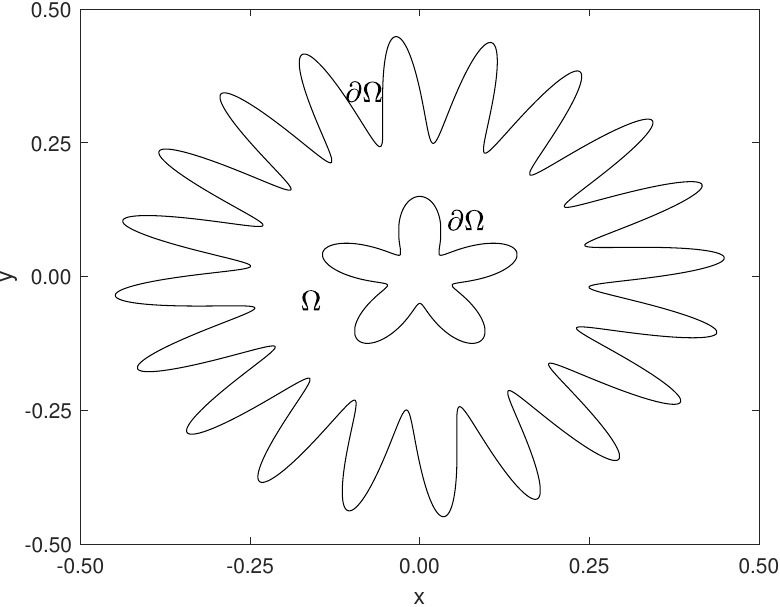}
		\end{subfigure}
		 \begin{subfigure}[b]{0.32\textwidth}
			 \includegraphics[width=6cm,height=6cm]{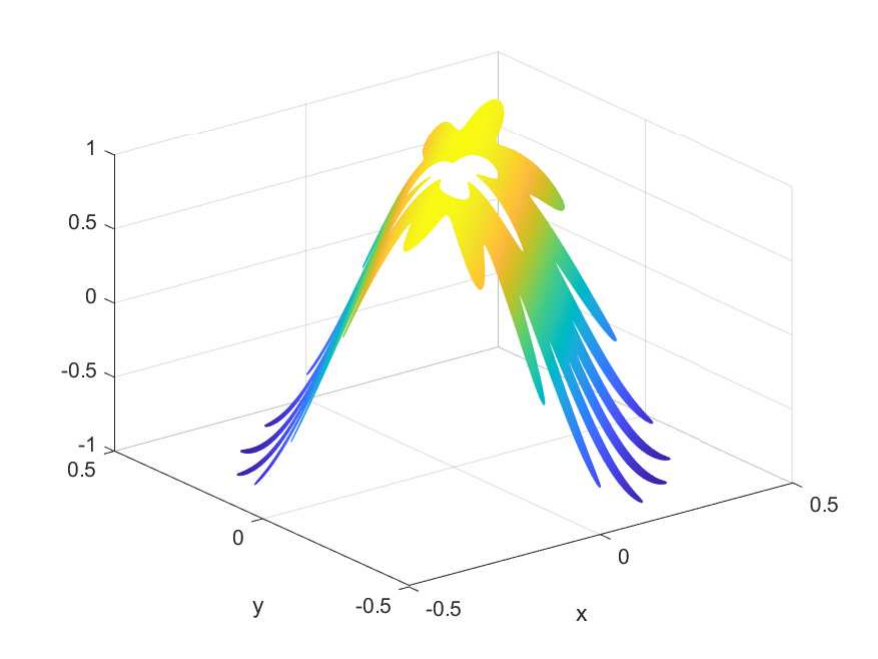}
		\end{subfigure}
		 \begin{subfigure}[b]{0.32\textwidth}
			 \includegraphics[width=6cm,height=6cm]{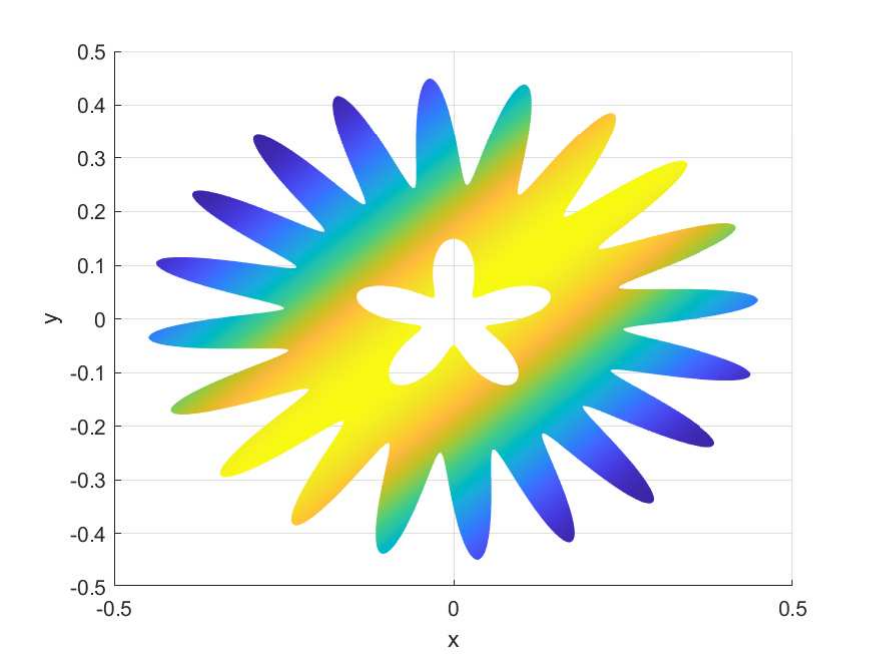}
		\end{subfigure}
		 \begin{subfigure}[b]{0.37\textwidth}
			 \includegraphics[width=6cm,height=6cm]{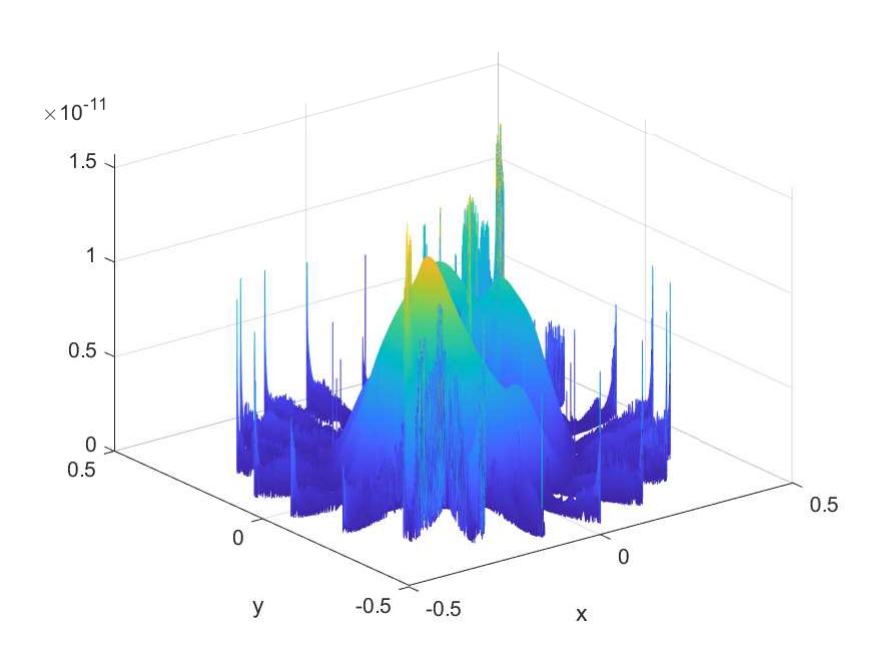}
		\end{subfigure}
		 \begin{subfigure}[b]{0.32\textwidth}
			 \includegraphics[width=6cm,height=6cm]{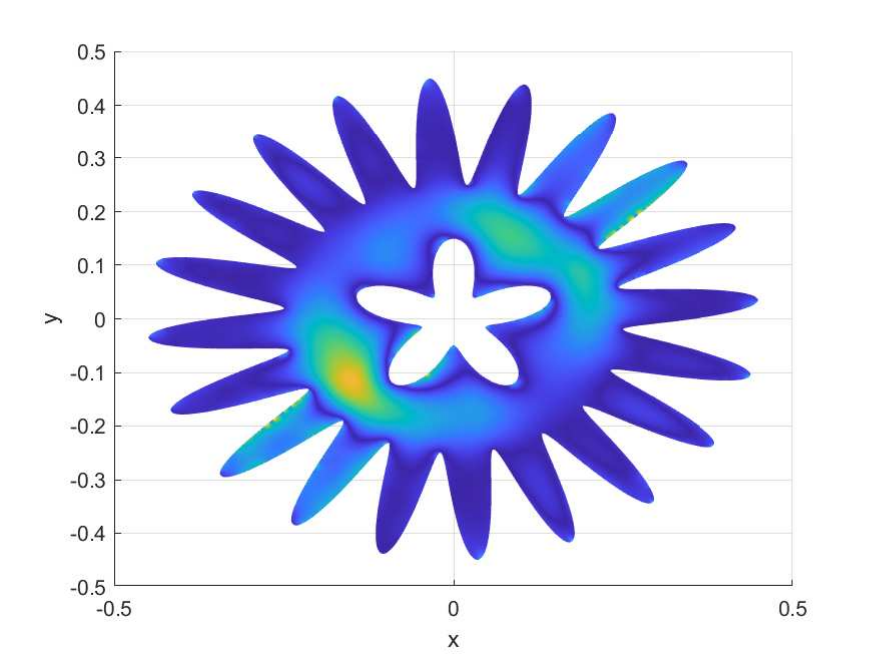}
		\end{subfigure}
		\caption
		{Performance in \cref{Example:2} of the proposed fourth-order compact FDM with $u=\cos(5(x-y))$. The domain $\Omega=\{(x(t), y(t)) : [0.1+0.05\sin(5t)]^2 < x(t)^2+y(t)^2 < [0.35+0.1\sin(20t)]^2, \ t\in [0,2\pi)\}$ that is enclosed by  5-leaf and 20-leaf boundary curves $\partial \Omega$ (left panel in the first row), $u_h$  on $\overline{\Omega}_h$ with $h=1/2^{12}$ (middle and right panels in the first row), and $|u_h-u|$   on $\overline{\Omega}_h$  with $h=1/2^{12}$ (left and right panels in the second row).}
		\label{Example:2:fig:2}
	\end{figure}

	\begin{example}\label{Example:2}
	\normalfont
	The functions  of the convection-diffusion-reaction equation \eqref{model:problem:2D} on a  curved domain $\Omega$ are given by
	\begin{align*}
		&u=\cos(k(x-y)),\quad k=5,10,\qquad \alpha= 4+\sin(4x)\cos(4y) ,\qquad  \\
		&\beta_1= \sin(x+2y), \qquad \beta_2= \cos(2x)\cos(3y), \qquad \kappa=\cos(3x+y), \\
		& \Omega=\{(x(t), y(t)) \ : \ [0.2+0.15\sin(5t)]^2 < x(t)^2+y(t)^2 < [0.3+0.16\sin(5t)]^2, \quad t\in [0,2\pi)\},\\
		& \Omega=\{(x(t), y(t)) \ : \ [0.1+0.05\sin(5t)]^2 < x(t)^2+y(t)^2 < [0.35+0.1\sin(20t)]^2, \quad t\in [0,2\pi)\},
	\end{align*}
	the source term $\phi$	and the Dirichlet boundary function $g$  are obtained by plugging above functions into \eqref{model:problem:2D}.  The numerical results are presented in \cref{Example:2:table:1,Example:2:table:2} and \cref{Example:2:fig:1,Example:2:fig:2}.	 
\end{example}

\begin{table}[htbp]
	\caption{Performance in \cref{Example:2} of the proposed fourth-order compact FDM, where $u=\cos(10(x-y))$ and $\Omega=\{(x(t), y(t)) \ : \ [0.2+0.15\sin(5t)]^2 < x(t)^2+y(t)^2 < [0.3+0.16\sin(5t)]^2, \ t\in [0,2\pi)\}$.}
	\centering
	{\renewcommand{\arraystretch}{1.0}
		\scalebox{1}{
			 \setlength{\tabcolsep}{5mm}{
				 \begin{tabular}{c|c|c|c|c}
					\hline
					$h$ &   $\frac{\|u_h-u\|_2}{\|u\|_2}$    &order &   	   $\|u_{h}-u\|_{\infty}$    &order \\
					\hline
$1/2^8$  &  2.7119E-04  &    &  1.1166E-02  &  \\
$1/2^9$  &  2.2578E-07  &  10.23  &  3.6442E-06  &  11.58\\
$1/2^{10}$  &  1.3317E-08  &  4.08  &  1.9743E-07  &  4.21\\
$1/2^{11}$  &  8.3414E-10  &  4.00  &  3.4156E-09  &  5.85\\
$1/2^{12}$  &  5.2193E-11  &  4.00  &  2.2667E-10  &  3.91\\
$1/2^{13}$  &  7.9523E-12  &  2.71  &  1.6872E-11  &  3.75\\
					\hline
	\end{tabular}}}}
	\label{Example:2:table:1}
\end{table}
		
	\begin{table}[htbp]
		\caption{Performance in \cref{Example:2} of the proposed fourth-order compact FDM, where $u=\cos(5(x-y))$ and $\Omega=\{(x(t), y(t)) \ : \ [0.1+0.05\sin(5t)]^2 < x(t)^2+y(t)^2 < [0.35+0.1\sin(20t)]^2, \ t\in [0,2\pi)\}$.}
		\centering
		 {\renewcommand{\arraystretch}{1.0}
			\scalebox{1}{
				 \setlength{\tabcolsep}{5mm}{
					 \begin{tabular}{c|c|c|c|c}
						\hline
						$h$ &   $\frac{\|u_h-u\|_2}{\|u\|_2}$    &order &   	   $\|u_{h}-u\|_{\infty}$    &order \\
						\hline
					$1/2^6$ &  1.9442E-02 &   &  2.1208E-01 &  \\
					$1/2^7$ &  1.0112E-03 &  4.26 &  2.4589E-02 &  3.11\\
					$1/2^8$ &  2.5933E-07 &  11.93 &  2.0397E-06 &  13.56\\
					$1/2^9$ &  1.5692E-08 &  4.05 &  5.1251E-08 &  5.31\\
					$1/2^{10}$ &  9.6204E-10 &  4.03 &  3.9733E-09 &  3.69\\
					$1/2^{11}$ &  6.3641E-11 &  3.92 &  2.4972E-10 &  3.99\\
					$1/2^{12}$ &  5.6048E-12 &  3.51 &  1.5665E-11 &  3.99\\	
						\hline
		\end{tabular}}}}
		\label{Example:2:table:2}
	\end{table}

\begin{figure}[htbp]
	\centering
	 \begin{subfigure}[b]{0.32\textwidth}
		 \includegraphics[width=5.7cm,height=5.7cm]{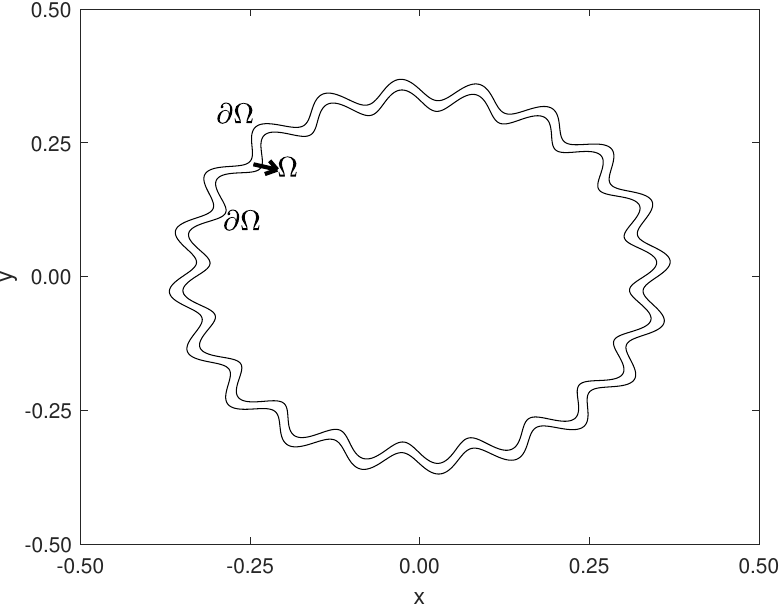}
	\end{subfigure}
	 \begin{subfigure}[b]{0.32\textwidth}
		 \includegraphics[width=6cm,height=6cm]{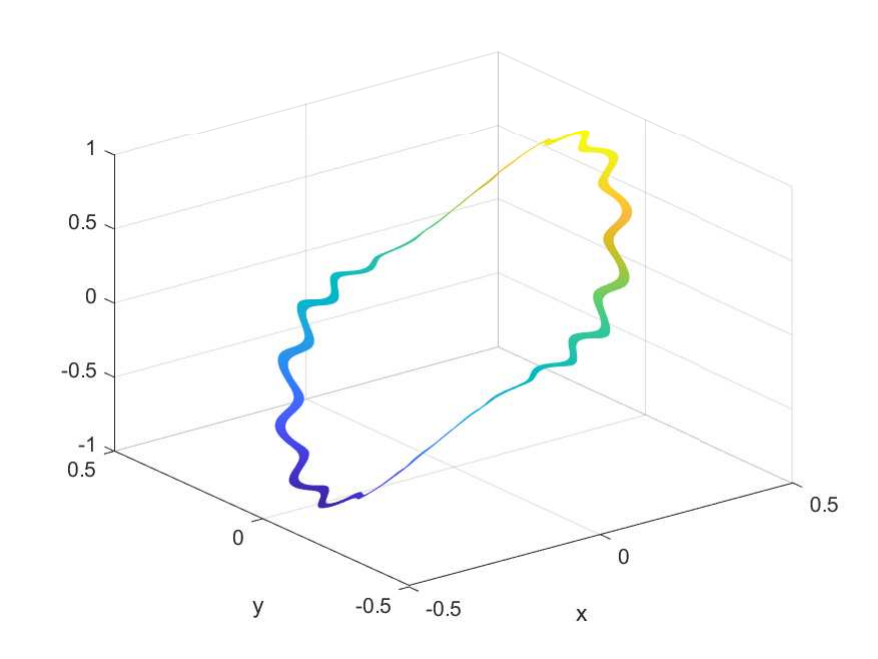}
	\end{subfigure}
	 \begin{subfigure}[b]{0.32\textwidth}
		 \includegraphics[width=6cm,height=6cm]{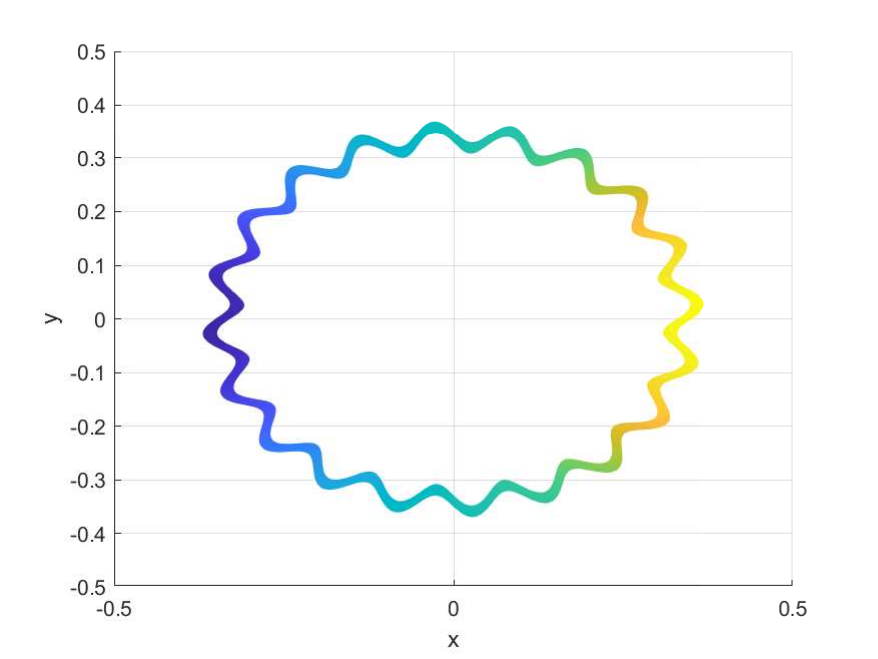}
	\end{subfigure}
	 \begin{subfigure}[b]{0.37\textwidth}
		 \includegraphics[width=6cm,height=6cm]{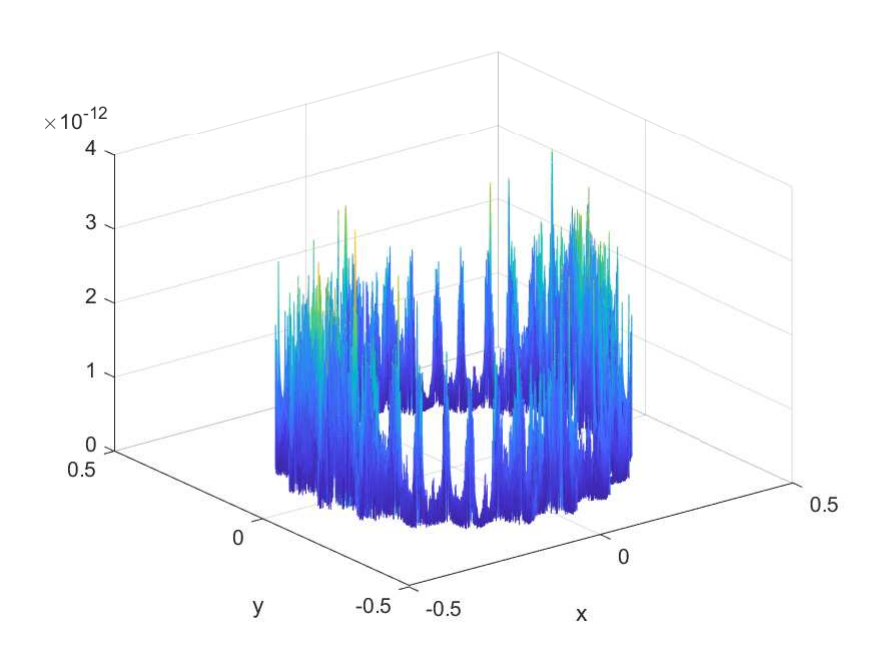}
	\end{subfigure}
	 \begin{subfigure}[b]{0.32\textwidth}
		 \includegraphics[width=6cm,height=6cm]{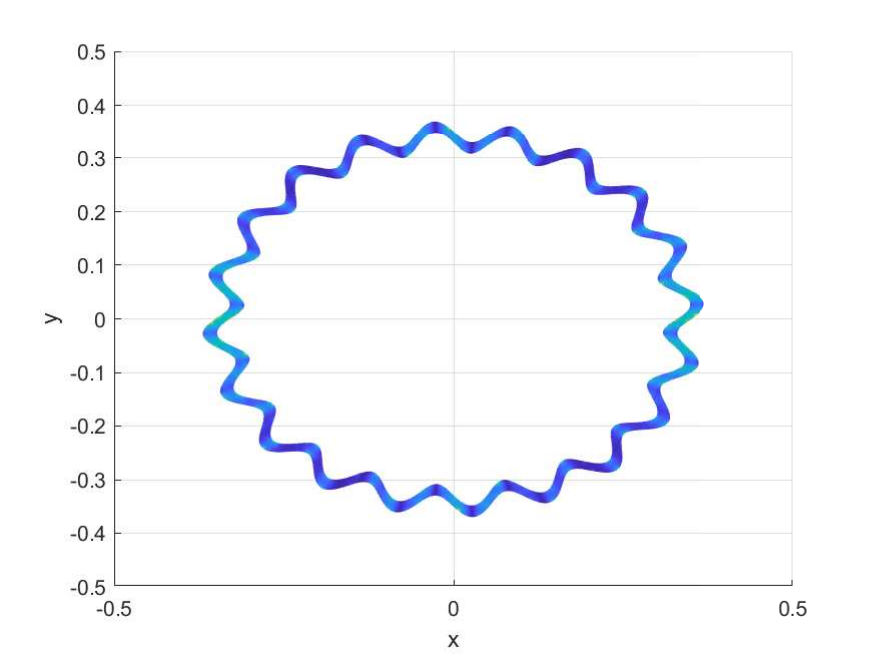}
	\end{subfigure}
	\caption
	{Performance in \cref{Example:3} of the proposed fourth-order compact FDM with  $u=\sin(4x)\cos(4y)$. The domain $\Omega=\{(x(t), y(t)) : [0.33+0.02\sin(20t)]^2 < x(t)^2+y(t)^2 < [0.35+0.02\sin(20t)]^2, \ t\in [0,2\pi)\}$ that is enclosed by  two 20-leaf boundary curves $\partial \Omega$ (left panel in the first row), $u_h$  on $\overline{\Omega}_h$ with $h=1/2^{12}$ (middle and right panels in the first row), and $|u_h-u|$   on $\overline{\Omega}_h$  with $h=1/2^{12}$ (left and right panels in the second row).}
	\label{Example:3:fig:1}
\end{figure}

\begin{figure}[htbp]
	\centering
	 \begin{subfigure}[b]{0.32\textwidth}
		 \includegraphics[width=5.7cm,height=5.7cm]{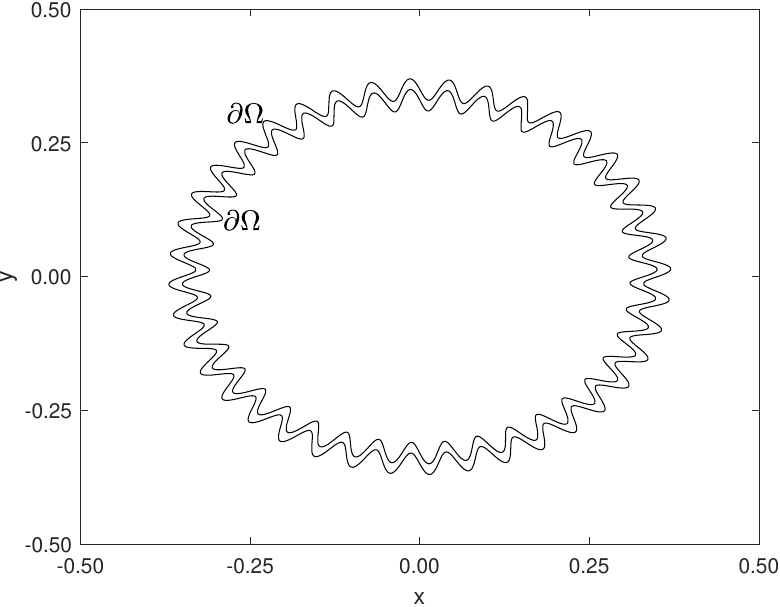}
	\end{subfigure}
	 \begin{subfigure}[b]{0.32\textwidth}
		 \includegraphics[width=6cm,height=6cm]{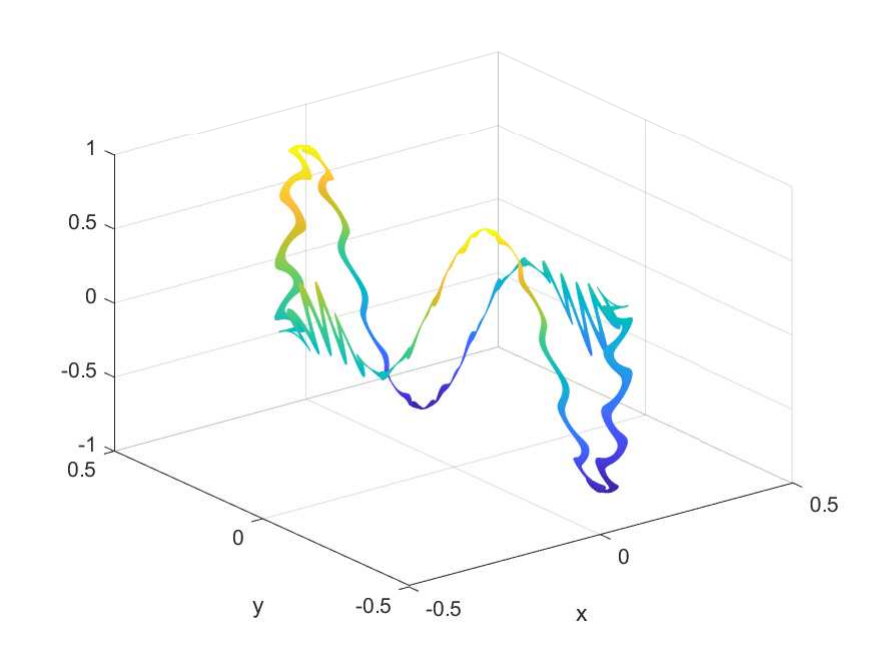}
	\end{subfigure}
	 \begin{subfigure}[b]{0.32\textwidth}
		 \includegraphics[width=6cm,height=6cm]{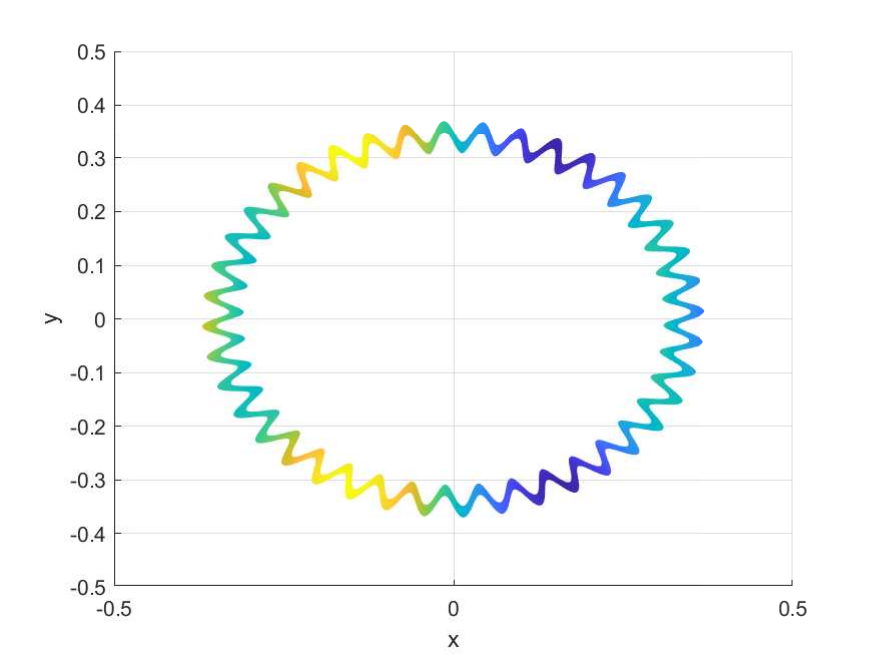}
	\end{subfigure}
	 \begin{subfigure}[b]{0.37\textwidth}
		 \includegraphics[width=6cm,height=6cm]{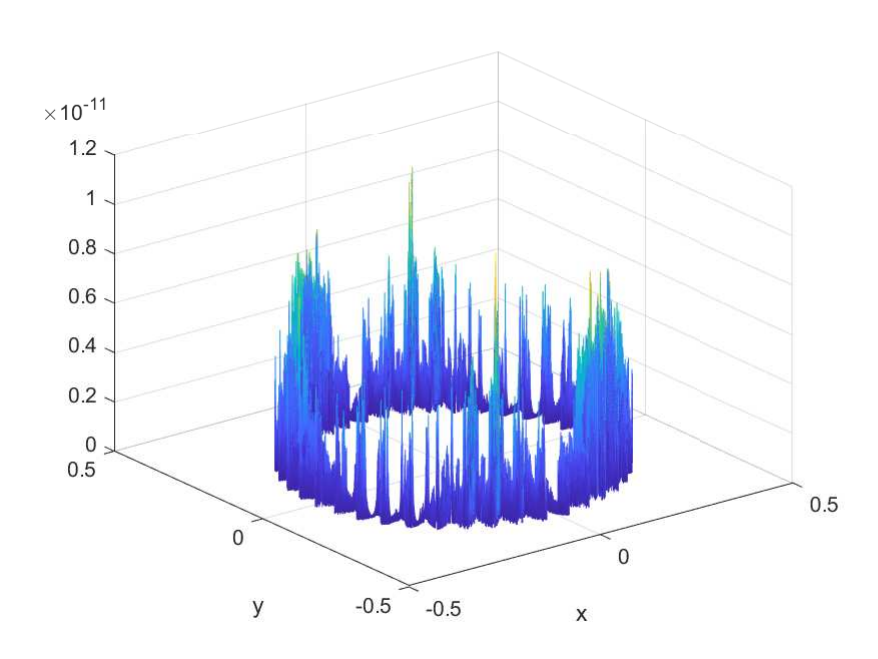}
	\end{subfigure}
	 \begin{subfigure}[b]{0.32\textwidth}
		 \includegraphics[width=6cm,height=6cm]{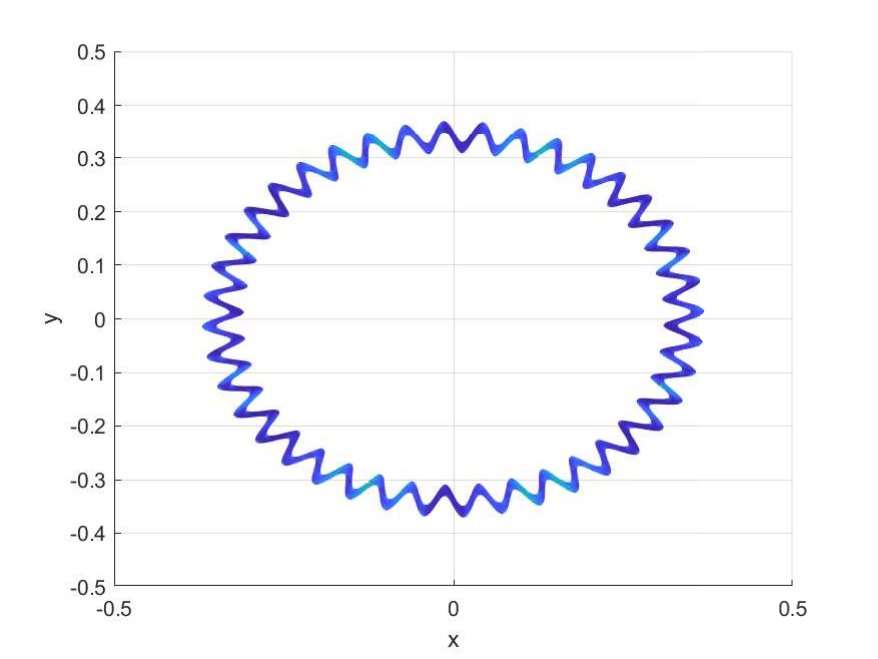}
	\end{subfigure}
	\caption
	{Performance in \cref{Example:3} of the proposed fourth-order compact FDM with $u=\sin(10x)\cos(10y)$. The domain $\Omega=\{(x(t), y(t)) : [0.33+0.02\sin(40t)]^2 < x(t)^2+y(t)^2 < [0.35+0.02\sin(40t)]^2, \ t\in [0,2\pi)\}$ that is enclosed by  two 40-leaf boundary curves $\partial \Omega$ (left panel in the first row), $u_h$  on $\overline{\Omega}_h$ with $h=1/2^{13}$ (middle and right panels in the first row), and $|u_h-u|$   on $\overline{\Omega}_h$  with $h=1/2^{13}$ (left and right panels in the second row).}
	\label{Example:3:fig:2}
\end{figure}

\begin{figure}[htbp]
	\centering
	 \begin{subfigure}[b]{0.32\textwidth}
		 \includegraphics[width=5.7cm,height=5.7cm]{curve6.pdf}
	\end{subfigure}
	 \begin{subfigure}[b]{0.32\textwidth}
		 \includegraphics[width=6cm,height=6cm]{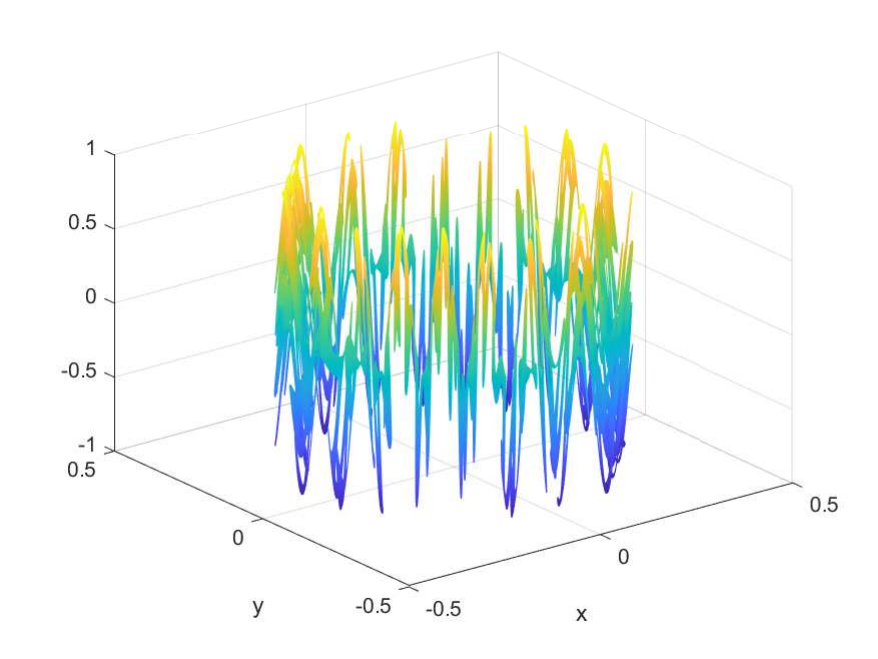}
	\end{subfigure}
	 \begin{subfigure}[b]{0.32\textwidth}
		 \includegraphics[width=6cm,height=6cm]{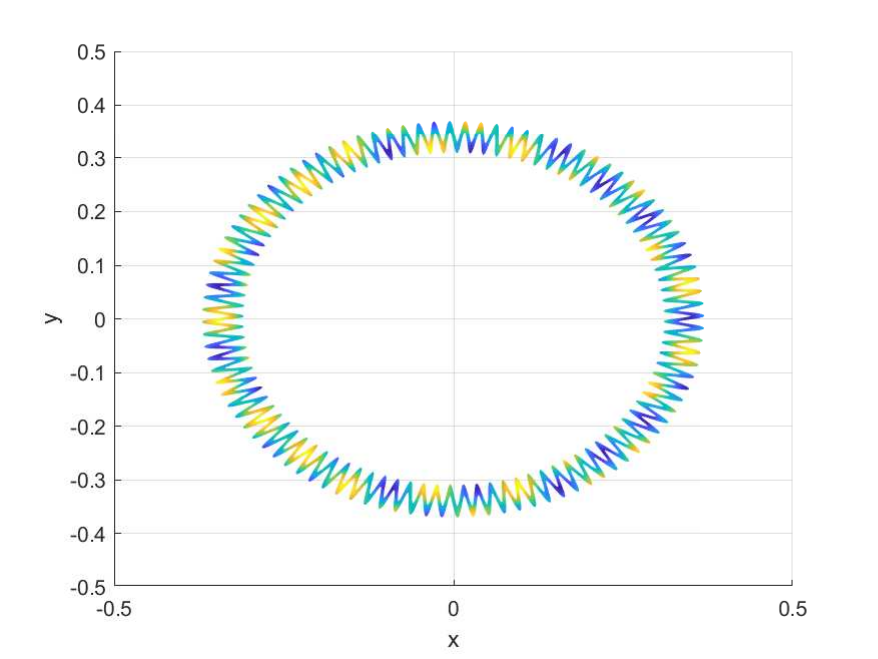}
	\end{subfigure}
	 \begin{subfigure}[b]{0.37\textwidth}
		 \includegraphics[width=6cm,height=6cm]{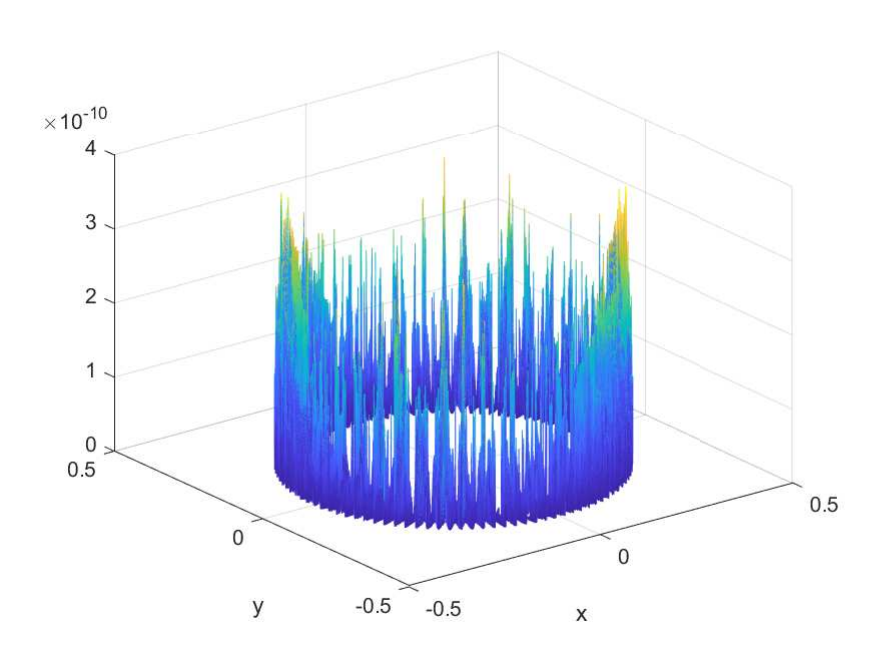}
	\end{subfigure}
	 \begin{subfigure}[b]{0.32\textwidth}
		 \includegraphics[width=6cm,height=6cm]{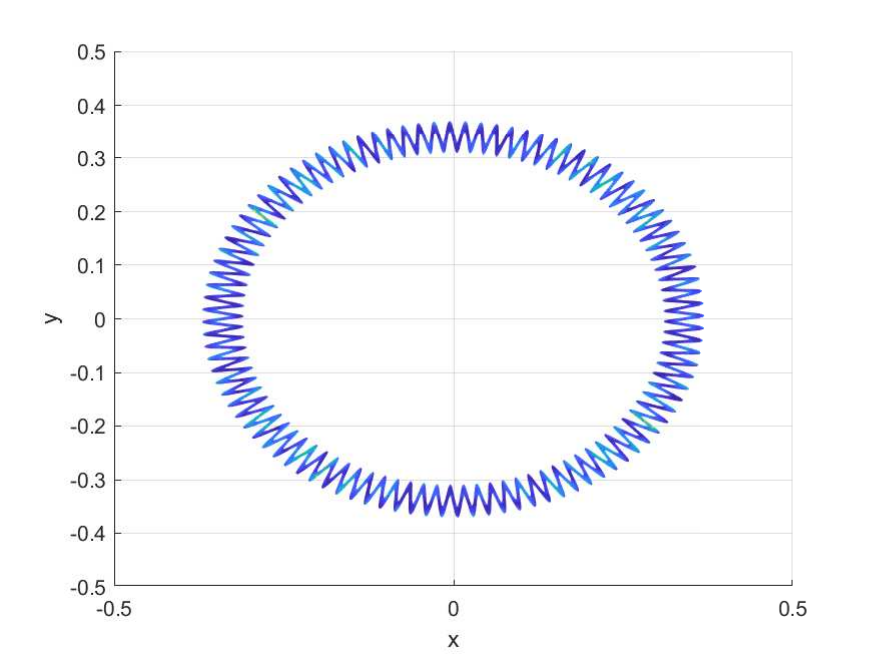}
	\end{subfigure}
	\caption
	{Performance in \cref{Example:3} of the proposed fourth-order compact FDM with $u=\sin(50x)\cos(50y)$. The domain  $\Omega=\{(x(t), y(t)) \ : \ [0.33+0.02\sin(100t)]^2 < x(t)^2+y(t)^2 < [0.35+0.02\sin(100t)]^2, \ t\in [0,2\pi)\}$ that is enclosed by  two 100-leaf boundary curves $\partial \Omega$ (left panel in the first row), $u_h$  on $\overline{\Omega}_h$ with $h=1/2^{14}$ (middle and right panels in the first row), and $|u_h-u|$   on $\overline{\Omega}_h$  with $h=1/2^{14}$ (left and right panels in the second row).}
	\label{Example:3:fig:3}
\end{figure}

\begin{figure}[htbp]
	\centering
	 \begin{subfigure}[b]{0.32\textwidth}
		 \includegraphics[width=5.7cm,height=5.7cm]{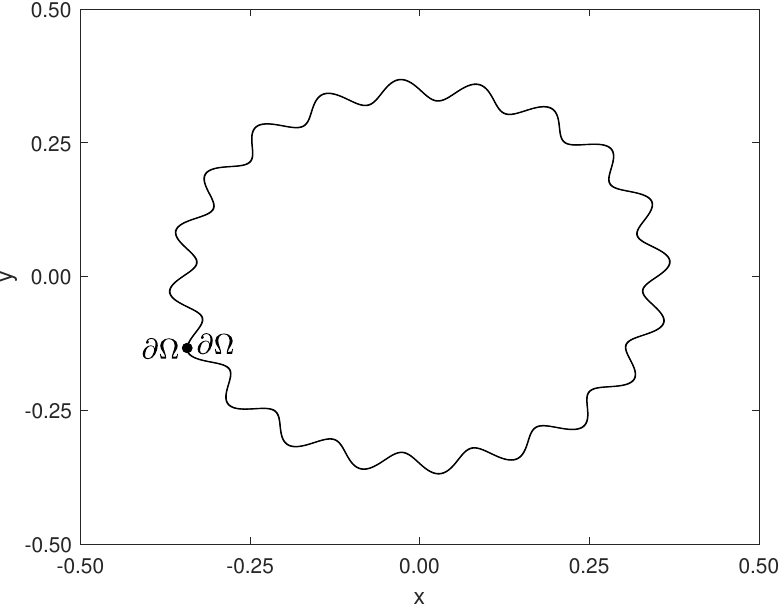}
	\end{subfigure}
	 \begin{subfigure}[b]{0.32\textwidth}
		 \includegraphics[width=5.6cm,height=5.6cm]{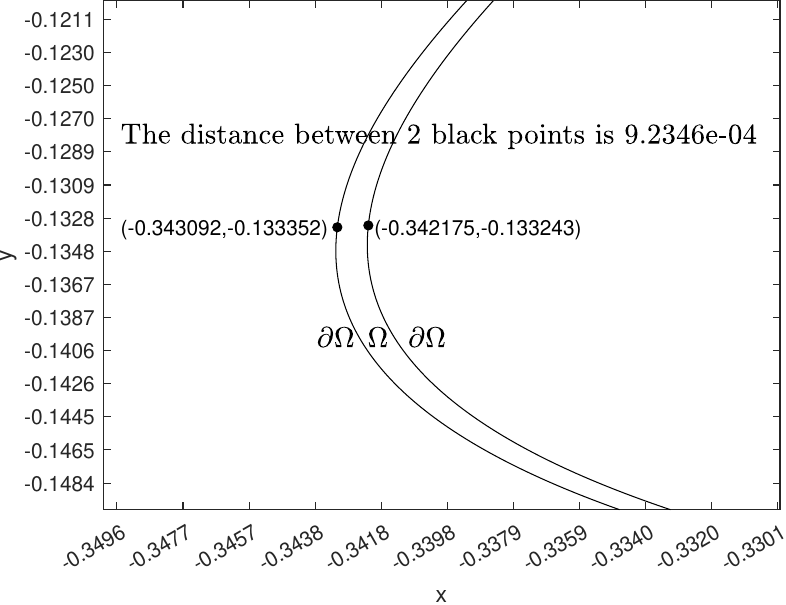}
	\end{subfigure}
	 \begin{subfigure}[b]{0.32\textwidth}
		 \includegraphics[width=6cm,height=6cm]{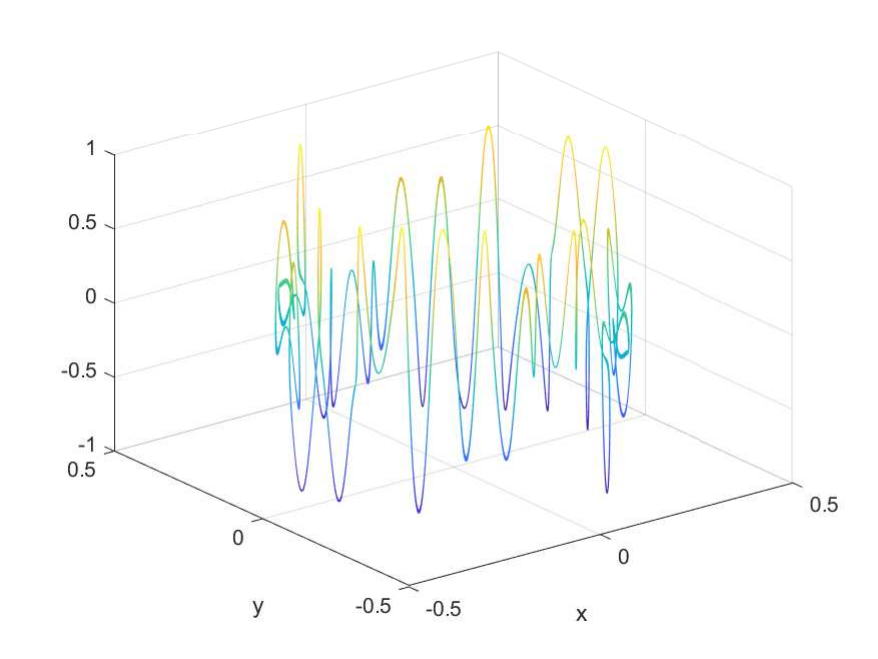}
	\end{subfigure}
	 \begin{subfigure}[b]{0.32\textwidth}
		 \includegraphics[width=6cm,height=6cm]{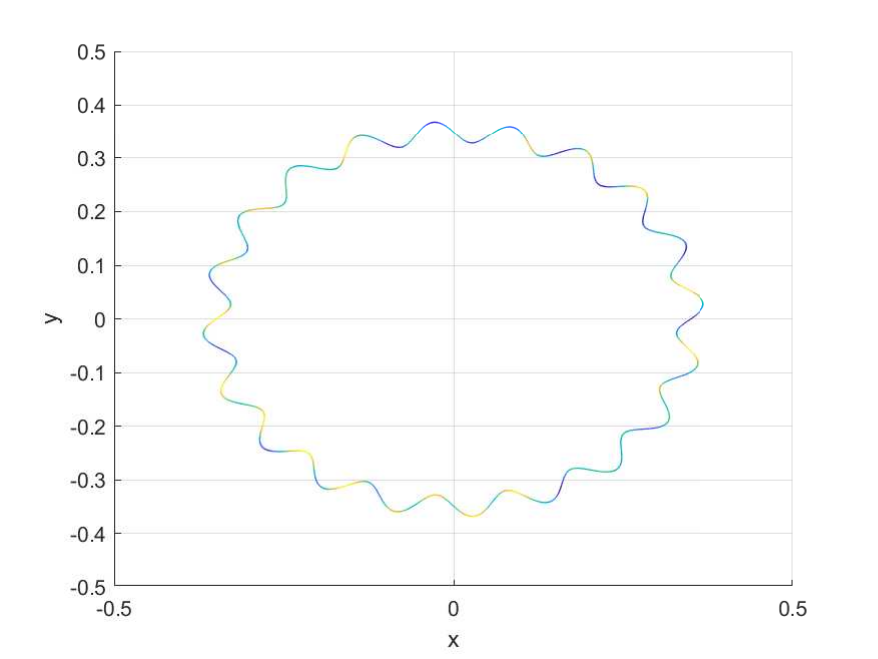}
	\end{subfigure}
	 \begin{subfigure}[b]{0.32\textwidth}
		 \includegraphics[width=6cm,height=6cm]{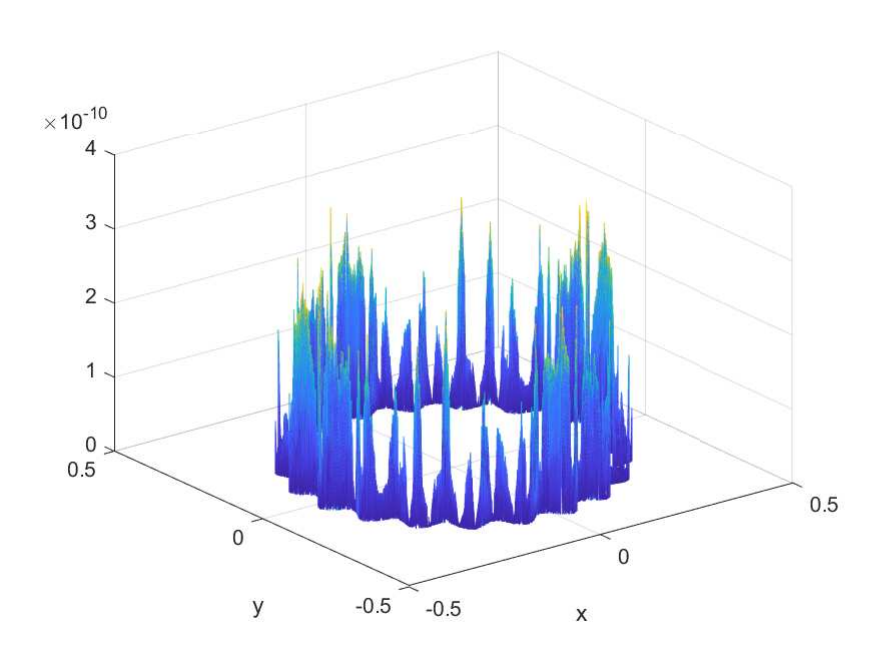}
	\end{subfigure}
	 \begin{subfigure}[b]{0.32\textwidth}
		 \includegraphics[width=6cm,height=6cm]{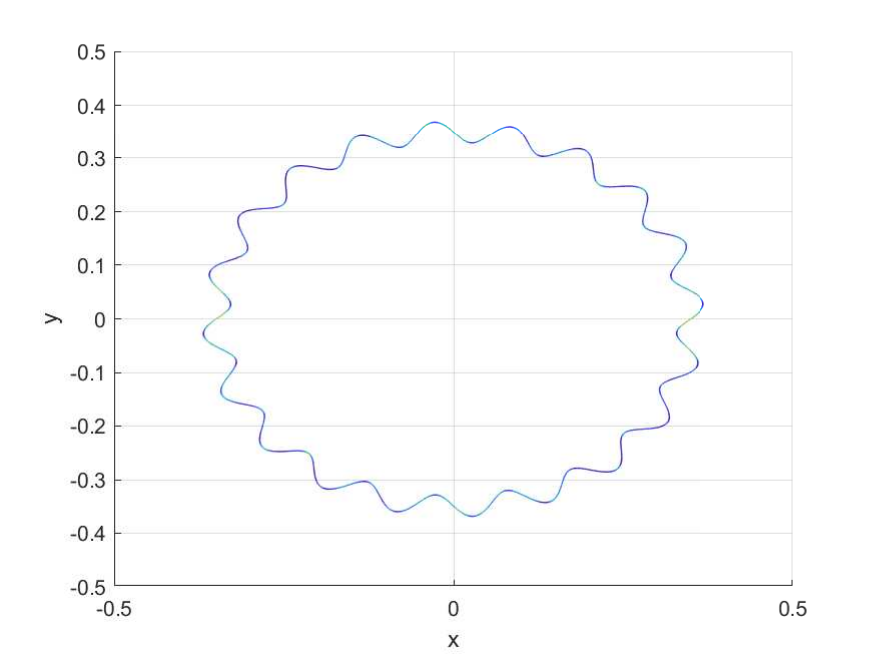}
	\end{subfigure}
	\caption
	{Performance in \cref{Example:3} of the proposed fourth-order compact FDM with $u=\sin(50x)\cos(50y)$. The domain  $\Omega=\{(x(t), y(t)) : [0.349+0.02\sin(20t)]^2 < x(t)^2+y(t)^2 < [0.35+0.02\sin(20t)]^2, \ t\in [0,2\pi)\}$ that is enclosed by  two nearly overlapping 20-leaf boundary curves $\partial \Omega$ (left and middle panels in the first row), $u_h$  on $\overline{\Omega}_h$ with $h=1/2^{14}$ (right panel in the first row and left panel in the second row), and $|u_h-u|$   on $\overline{\Omega}_h$  with $h=1/2^{14}$ (middle and right panels in the second row).}
	\label{Example:3:fig:4}
\end{figure}

	\begin{example}\label{Example:3}
		\normalfont
		The functions  of the convection-diffusion-reaction equation \eqref{model:problem:2D} on a curved domain $\Omega$ are given by
		\begin{align*}
		&u=\sin(k x)\cos(k y),\qquad k=4,10,50, \qquad \alpha= \exp(x+y) ,\qquad  \\
	&\beta_1= \exp(x-y), \qquad \beta_2= \cos(x)\cos(y), \qquad \kappa=\exp(x+y), \\
			& \Omega=\{(x(t), y(t)) \ : \ [0.33+0.02\sin(20t)]^2 < x(t)^2+y(t)^2 < [0.35+0.02\sin(20t)]^2, \quad t\in [0,2\pi)\},\\
			& \Omega=\{(x(t), y(t)) \ : \ [0.33+0.02\sin(40t)]^2 < x(t)^2+y(t)^2 < [0.35+0.02\sin(40t)]^2, \quad t\in [0,2\pi)\},\\
			& \Omega=\{(x(t), y(t)) \ : \ [0.33+0.02\sin(100t)]^2 < x(t)^2+y(t)^2 < [0.35+0.02\sin(100t)]^2, \quad t\in [0,2\pi)\},\\
			& \Omega=\{(x(t), y(t)) \ : \ [0.349+0.02\sin(20t)]^2 < x(t)^2+y(t)^2 < [0.35+0.02\sin(20t)]^2, \quad t\in [0,2\pi)\},
		\end{align*}
		the source term $\phi$	and the Dirichlet boundary function $g$  are obtained by plugging above functions into \eqref{model:problem:2D}. The numerical results are presented in \cref{Example:3:table:1,Example:3:table:2,Example:3:table:3,Example:3:table:4} and \cref{Example:3:fig:1,Example:3:fig:2,Example:3:fig:3,Example:3:fig:4}.	 
	\end{example}

	\begin{table}[htbp]
		\caption{Performance in \cref{Example:3} of the proposed fourth-order compact FDM, where $u=\sin(4x)\cos(4y)$ and $\Omega=\{(x(t), y(t))  :  [0.33+0.02\sin(20t)]^2 < x(t)^2+y(t)^2 < [0.35+0.02\sin(20t)]^2, \ t\in [0,2\pi)\}$.}
		\centering
		 {\renewcommand{\arraystretch}{1.0}
			\scalebox{1}{
				 \setlength{\tabcolsep}{5mm}{
					 \begin{tabular}{c|c|c|c|c}
						\hline
						$h$ &   $\frac{\|u_h-u\|_2}{\|u\|_2}$    &order &   	   $\|u_{h}-u\|_{\infty}$    &order \\
						\hline
$1/2^6$ &  1.1561E-02 &   &  4.7770E-02 &  \\
$1/2^7$ &  6.3478E-07 &  14.15 &  2.7195E-06 &  14.10\\
$1/2^8$ &  8.3083E-08 &  2.93 &  1.7030E-07 &  4.00\\
$1/2^9$ &  5.6609E-09 &  3.88 &  1.1265E-08 &  3.92\\
$1/2^{10}$ &  3.9934E-10 &  3.83 &  8.8971E-10 &  3.66\\
$1/2^{11}$ &  2.4814E-11 &  4.01 &  5.6864E-11 &  3.97\\
$1/2^{12}$ &  1.4931E-12 &  4.05 &  3.7550E-12 &  3.92\\
						\hline
		\end{tabular}}}}
		\label{Example:3:table:1}
	\end{table}
		
	\begin{table}[htbp]
		\caption{Performance in \cref{Example:3} of the proposed fourth-order compact FDM, where $u=\sin(10x)\cos(10y)$ and $\Omega=\{(x(t), y(t)) : [0.33+0.02\sin(40t)]^2 < x(t)^2+y(t)^2 < [0.35+0.02\sin(40t)]^2, \ t\in [0,2\pi)\}$.}
		\centering
		 {\renewcommand{\arraystretch}{1.0}
			\scalebox{1}{
				 \setlength{\tabcolsep}{5mm}{
					 \begin{tabular}{c|c|c|c|c}
						\hline
						$h$ &   $\frac{\|u_h-u\|_2}{\|u\|_2}$    &order &   	   $\|u_{h}-u\|_{\infty}$    &order \\
						\hline
$1/2^7$ &  3.0546E-02 &   &  3.1862E-01 &  \\
$1/2^8$ &  2.7615E-06 &  13.43 &  1.6867E-05 &  14.21\\
$1/2^9$ &  1.7321E-07 &  3.99 &  5.0266E-07 &  5.07\\
$1/2^{10}$ &  1.1805E-08 &  3.87 &  3.1117E-08 &  4.01\\
$1/2^{11}$ &  7.8057E-10 &  3.92 &  2.2230E-09 &  3.81\\
$1/2^{12}$ &  5.0312E-11 &  3.96 &  1.4953E-10 &  3.89\\
$1/2^{13}$ &  2.9737E-12 &  4.08 &  1.0823E-11 &  3.79\\
						\hline
		\end{tabular}}}}
		\label{Example:3:table:2}
	\end{table}

		\begin{table}[htbp]
		\caption{Performance in \cref{Example:3} of the proposed fourth-order compact FDM, where $u=\sin(50x)\cos(50y)$ and $\Omega=\{(x(t), y(t))  :  [0.33+0.02\sin(100t)]^2 < x(t)^2+y(t)^2 < [0.35+0.02\sin(100t)]^2, \ t\in [0,2\pi)\}$.}
		\centering
		 {\renewcommand{\arraystretch}{1.0}
			\scalebox{1}{
				 \setlength{\tabcolsep}{5mm}{
					 \begin{tabular}{c|c|c|c|c}
						\hline
						$h$ &   $\frac{\|u_h-u\|_2}{\|u\|_2}$    &order &   	   $\|u_{h}-u\|_{\infty}$    &order \\
						\hline
$1/2^8$  &  5.7873E-01  &    &  6.4646E+00  &  \\
$1/2^9$  &  1.5027E-01  &  1.95  &  3.2825E+00  &  0.98\\
$1/2^{10}$  &  4.2131E-03  &  5.16  &  1.3816E-01  &  4.57\\
$1/2^{11}$  &  5.3721E-07  &  12.94  &  3.2870E-06  &  15.36\\
$1/2^{12}$  &  3.4059E-08  &  3.98  &  8.8488E-08  &  5.22\\
$1/2^{13}$  &  2.1544E-09  &  3.98  &  5.6768E-09  &  3.96\\
$1/2^{14}$  &  1.3652E-10  &  3.98  &  3.7503E-10  &  3.92\\
						\hline
		\end{tabular}}}}
		\label{Example:3:table:3}
	\end{table}

		\begin{table}[htbp]
		\caption{Performance in \cref{Example:3} of the proposed fourth-order compact FDM, where $u=\sin(50x)\cos(50y)$ and $\Omega=\{(x(t), y(t)) : [0.349+0.02\sin(20t)]^2 < x(t)^2+y(t)^2 < [0.35+0.02\sin(20t)]^2, \  t\in [0,2\pi)\}$.}
		\centering
		 {\renewcommand{\arraystretch}{1.0}
			\scalebox{1}{
				 \setlength{\tabcolsep}{5mm}{
					 \begin{tabular}{c|c|c|c|c}
						\hline
						$h$ &   $\frac{\|u_h-u\|_2}{\|u\|_2}$    &order &   	   $\|u_{h}-u\|_{\infty}$    &order \\
						\hline
$1/2^{10}$  &  1.3975E-02  &    &  1.8476E-01  &  \\
$1/2^{11}$  &  1.8983E-03  &  2.88  &  1.5374E-02  &  3.59\\
$1/2^{12}$  &  3.0217E-08  &  15.94  &  7.9659E-08  &  17.56\\
$1/2^{13}$  &  2.3576E-09  &  3.68  &  5.5330E-09  &  3.85\\
$1/2^{14}$  &  1.5499E-10  &  3.93  &  3.2938E-10  &  4.07\\
						\hline
		\end{tabular}}}}
		\label{Example:3:table:4}
	\end{table}

	 \section{Contribution}\label{sec:contribu}
	
In this paper, we consider the convection-diffusion-reaction equation with smooth variable convection, diffusion, reaction functions and the Dirichlet boundary condition on a smooth curved domain.  We propose the fourth-order compact FDM to approximate the solution at regular and irregular  stencil center points using the uniform Cartesian grid. For the thin domain with the sharply varying boundary curve, there are many configurations of the compact FDM at the irregular  stencil center. We establish a linear system with an at most $6 \times 24$ matrix to compute the left-hand side of the fourth-order compact FDM to cover all configurations. Furthermore, each entry of the matrix is expressed in the explicit expression for any irregular  stencil centers.  Once the left-hand side is computed, the right-hand side of the stencil is calculated immediately by the given explicit expression. \\
Strengths of our FDM:
\begin{itemize}
	\item All formulas are written in explicit forms such that our FDM is simple to be implemented to help readers reproduce our results straightforwardly.
	\item From numerical examples, our FDM performs well for the boundary curve with the large curvature (e.g. \cref{Example:1:fig,Example:2:fig:1} and \cref{Example:1:table,Example:2:table:1}) and the thin domain that cannot be distinguished by the naked eye (e.g. \cref{Example:3:fig:4} and \cref{Example:3:table:4}). Our FDM yields the stable convergence rate and the accurate solution for the high-frequency solution and the highly-oscillatory curve (e.g. \cref{Example:3:fig:3} and \cref{Example:3:table:3}). Furthermore, all errors are uniformly distributed around the boundary, which demonstrates the stability and robustness of our FDM.
	\item From results in all tables, we observe that the convergence rate is stabilized at 4 when $h$ is reasonably small, which numerically indicates that the convergence of our proposed FDM at the irregular  stencil center seems to be independent of the distance between the centered grid point and the boundary curve, or this distance has the limited effect on the accuracy of our FDM.
\end{itemize}
Key differences from existing methods:
\begin{itemize}
	\item In contrast to the ghost cell method used in \cite{Clain2021,Clain2024}, we do not need to use approximated values of $u$ at grid points outside $\Omega$ by the interpolation.
	\item Compared with the immersed interface method (IIM) in \cite {Ito2005,LiIto06}, we do not utilize the coordinate transformation and solve the optimization problem to obtain the stencil. We use the transformation in Cartesian coordinate directly to achieve the efficient implementation, and solve a small linear system with the explicit expression to compute stencil to improve the accuracy. Furthermore, \cite {Ito2005,LiIto06} employ 1 to 4 points at the boundary in the stencil at the irregular  stencil center, and we use one point at the boundary to derive the fourth-order compact FDM.
	\item  Distinct from using partial derivatives in the complex field to derive the FDM in \cite{HanSim2025}, we only use the real partial derivatives to build the simple stencil and enhance the efficiency.
	\item Different from the augmented matched interface and boundary (AMIB) method in \cite{Li2023,Li2025,Ren2022},
	we construct the compact scheme which generates a sparse matrix with nine non-zero bands.
\end{itemize}
The potential of our method:
\begin{itemize}
	\item Using techniques in \cite{FengTrenchea2026,Feng3D2026}, it is direct to extend the proposed fourth-order compact FDM to solve
	nonlinear steady and unsteady convection-diffusion-reaction equations on irregular domains in 2D and 3D.
	\item Due to the simple and efficient implementation of our method, it is also straightforward to include more points in the stencil to attain the higher-order accuracy for more challenging PDEs (e.g, the 25-point FDM for the Stokes problem in \cite{FengHanNeilan2026}). Furthermore, it is also feasible to extend our FDM to solve the parabolic equation with the moving boundary/interface and the free boundary problem by the flexibility and efficiency of our derivation.
	\item We can derive the high-order FDM for \eqref{model:problem:2D} with Neumann and Robin boundary conditions by replacing \eqref{Dirichlet} by $\rho u+ \zeta \frac{\partial u}{\partial \vec{n}}=g$ and repeating  \eqref{f00}--\eqref{Fij:explicit} similarly.
	\item The proposed FDM with more points in the stencil only results in more matrices $A_{r,\ell}$ in the linear system \eqref{sub:linear:system}. So the
	high-order FDM with the adaptive mesh can be derived easily by adding  more $A_{r,\ell}$ in \eqref{sub:linear:system} for the irregular  stencil center near the sharply varying boundary curve.
\end{itemize}

	\section{Declarations}
	\noindent \textbf{Conflict of interest:} The authors declare that they have no conflict of interest.\\
	\noindent \textbf{Data availability:} Data will be made available on reasonable request.
	
	\vspace{0.3cm}
	\noindent\textbf{Acknowledgment}
	
	Qiwei Feng is partially supported by   the Mathematics Department, King Fahd University of Petroleum and Minerals (KFUPM), Dhahran, Saudi Arabia.
	
	Bin Han is partially supported by Natural Sciences and Engineering Research Council (NSERC) of Canada under grants RGPIN-2024-04991.

The main idea of this manuscript was developed in 2021. Owing to various reasons, the work was delayed. Regrettably, the manuscript was not completed before the third author Peter Minev left us.

\end{document}